\pdfoutput=1
\documentclass[11pt]{amsart}
\usepackage[all]{background}
\backgroundsetup{contents={}}
\usepackage[utf8]{inputenc}
\usepackage{changepage}
\usepackage[margin=1in]{geometry}  
\usepackage{graphicx}              
\usepackage{color}              
\usepackage{float}
\usepackage{amsmath}               
\usepackage{amssymb}               
\usepackage{amsfonts}              
\usepackage{amsthm}                
\usepackage{tikz}                  
\usepackage[pdfencoding=auto, psdextra]{hyperref}                  
\usepackage{microtype}                  
\usepackage{xparse}

\usepackage{algorithm}              
\usepackage{algorithmicx} 
\usepackage{algpseudocode} 

\usepackage{microtype}
\usepackage{verbatim}
\usepackage{tikz}
\usetikzlibrary{arrows, decorations.markings}
\usetikzlibrary{calc} 
\usepackage{longtable}                  
\usepackage{lineno}                  
\usetikzlibrary{matrix}
\usepackage{enumitem} 
\setlist[enumerate]{topsep=0pt,itemsep=-1ex,partopsep=1ex,parsep=1ex}

\theoremstyle{plain}
\newtheorem{theorem}{Theorem}[section]
\newtheorem{lemma}[theorem]{Lemma}
\newtheorem{proposition}[theorem]{Proposition}
\newtheorem{corollary}[theorem]{Corollary}

\newtheorem{definition}[theorem]{Definition}
\newtheorem{example}[theorem]{Example}

\newtheorem*{theorem*}{Theorem}
\newtheorem*{lemma*}{Lemma}
\newtheorem*{proposition*}{Proposition}
\newtheorem*{corollary*}{Corollary}
\newtheorem*{conjecture*}{Conjecture}
\newtheorem*{definition*}{Definition}
\newtheorem*{example*}{Example}
\newtheorem*{remark*}{Remark}
\newtheorem*{question*}{Question}

\newcommand{\lieGroup}{\mathsf} 
 
\newcommand{\mathSet}{\mathrm} 

\ExplSyntaxOn
\NewDocumentCommand{\captureshell}{ov}
 {
  \IfNoValueTF { #1 }
   {
    \sdaau_captureshell:Nn \l__sdaau_captureshell_out_tl { #2 }
    \tl_use:N \l__sdaau_captureshell_out_tl
   }
   {
    \sdaau_captureshell:Nn #1 { #2 }
   }
 }

\tl_new:N \l__sdaau_captureshell_out_tl

\cs_new_protected:Nn \sdaau_captureshell:Nn
 {
  \tl_set_from_file:Nnn #1 { } { |"#2" }
 }
\ExplSyntaxOff

\usetikzlibrary{patterns}
\title[Disconnecting the moduli space of $\lieGroup{G}_2$-metrics]{Disconnecting the moduli space of $\lieGroup{G}_2$-metrics via $\lieGroup{U}(4)$-coboundary defects}
\author{Dominic Wallis}\thanks{D. Wallis, University of Bath, \href{mailto:dw580@bath.ac.uk}{dw580@bath.ac.uk}}
\hypersetup{
    pdftitle = {Disconnecting the moduli space of G2 metrics via U(4)-coboundary defects},
    pdfauthor = {Dominic Wallis},
    pdfsubject = {},
    pdfkeywords = {AMS MSC. Primary: 53C10, 57R15; Secondary: 53C25, 53C27, 53D15, 57R90.}
}

\begin{document} 

\begin{abstract}
    We exhibit examples of closed Riemannian 7-manifolds with holonomy $\lieGroup{G}_2$ such that the underlying manifolds are 
    diffeomorphic but whose associated $\lieGroup{G}_2$-structures are not homotopic. 
    This is achieved by defining two invariants of certain $\lieGroup{U}(3)$-structures. 
    We show that these agree with the invariants of $\lieGroup{G}_2$-structures defined by Crowley and Nordstr\"om. 
    We construct a suitable coboundary for $\lieGroup{G}_2$ manifolds obtained via the Twisted Connected Sum method 
    that allows the invariants to be computed in terms of the input data of the construction.  
    We find explicit examples where the invariants detect different connected components of the moduli of $\lieGroup{G}_2$-metrics.
\end{abstract}

\maketitle

\section{Introduction} 

We exhibit examples of closed $7$-manifolds for which the moduli space of metrics with holonomy precisely $\lieGroup{G}_2$ is disconnected. 
By $\lieGroup{G}_2$-metric, we mean a metric with holonomy precisely $\lieGroup{G}_2$, which together form a moduli space. 
A $\lieGroup{G}_2$-metric induces a reduction of the structure group to a $\lieGroup{G}_2$-structure. 
Two $G$-structures are homotopic if there exists a path between them through $G$-structures.
In particular, if two $\lieGroup{G}_2$-metrics belong to the same connected component of the moduli space, 
then their associated $\lieGroup{G}_2$-structures are homotopic up to diffeomorphism. 
Note that the converse is not necessarily true.  

We begin by defining two invariants of certain $\lieGroup{U}(3)$-structures on 7-manifolds denoted $\nu$ and $\xi$. 
The invariants are constant on homotopy classes of $\lieGroup{U}(3)$-structures. 
The definitions are in terms of a $\lieGroup{U}(4)$-coboundary (ie a coboundary with $\lieGroup{U}(4)$-structure that agrees with the $\lieGroup{U}(3)$-structure on the boundary)
and are understood as `boundary defects' of almost complex 8-manifolds.
We show that the definitions are well defined in the sense that a suitable coboundary exists 
and that the definitions are independent of choices made.
This boundary defect approach was most notably used by Milnor \cite{milnor1956exotic} when he employed his $\lambda$-invariant to prove the existence of exotic 7-spheres.

Crowley and Nordstr\"om \cite{crowley2015new} define invariants of $\lieGroup{G}_2$-structures on 7-manifolds. 
By noting a relationship between $\lieGroup{G}_2$-structures and $\lieGroup{U}(3)$-structures on 7-manifolds, 
we show that our invariants are equivalent to those of \cite{crowley2015new} in the appropriate context. 
This gives an alternative method of computing the invariants that had previously been intractable in certain examples of interest.  

The second half of the paper is focused on finding examples for which our invariants demonstrate the claimed phenomenon: separating the moduli of $\lieGroup{G}_2$-metrics. 
There are essentially two known constructions of $\lieGroup{G}_2$-manifolds: Joyce's method of desingularizing quotient tori \cite{joyce1996compact}; and the twisted connected sum (TCS) method of Kovalev \cite{kovalev2000twisted}, 
extended by Corti et al \cite{corti2015}. 
We focus on the latter, which constructs $\lieGroup{G}_2$-manifolds from pairs of certain complex threefolds called building blocks.
Of the two invariants defined in \cite{crowley2015new}, one is unable to distinguish TCS manifolds (Proposition \ref{prop_nu_is_24}),
while the other has not previously been computed.
Examples of closed $7$-manifolds for which the moduli space of $\lieGroup{G}_2$-metrics is disconnected are known \cite{crowley2015analytic},
but these examples are obtained by extending the TCS method with `extra twisting'. 

We describe a TCS $\lieGroup{U}(4)$-coboundary that allows us to compute our invariants for TCS manifolds.
We reformulate $\nu$ and $\xi$ in terms of data of building blocks and some arrangement of lattices (Definition \ref{def_configuration}) that essentially determines the TCS.
From this point $\nu$ and $\xi$ can be computed for any TCS from cohomological data of building blocks, and this lattice arrangement.

One way to obtain building blocks is via weak Fano threefolds. 
Thanks to advances in our understanding of weak Fanos (such as \cite{arap2017existence} and references therein) much of the data 
required to compute relevant invariants can be read from the literature with little extra work. 
By classification results of 2-connected 7-manifolds in terms of topological invariants \cite{crowley2014classification}, 
we can determine the diffeomorphism class of a TCS. 

We conclude with the first known examples of pairs of TCS manifolds which have underlying diffeomorphic smooth structures, but nonhomotopic $\lieGroup{G}_2$-structures. 
Thus the moduli space of $\lieGroup{G}_2$-metrics over these smooth manifolds is disconnected.

\subsection{The invariants} 

Let $M$ be a smooth closed 7-manifold, with a $\lieGroup{U}(3)$-structure. 
A $\lieGroup{U}(4)$-coboundary $W$ of $M$, is a compact 8-manifold with $\partial W \cong M$ and with a $\lieGroup{U}(4)$-structure,   
such that on restriction to the boundary it agrees with the $\lieGroup{U}(3)$-structure on $M$. 
We say that a coboundary $W$ of $M$ is over $H^4$ if the map $H^4(W) \rightarrow H^4(M)$ is surjective.
Throughout cohomology is assumed to be with integer coefficients if not otherwise stated. 

\begin{definition}\label{def_nu_u3}
    Let $M$ be a closed 7-manifold with $\lieGroup{U}(3) $-structure $(v, g, \omega)$.
    Assume that the Chern classes $c_1(M) \in H^2(M)$ and $c_3 (M) \in H^6(M)$ both vanish.
    Let $W$ be a $\lieGroup{U}(4)$-coboundary of $M$ over $H^4$. Then 
    \begin{equation}\begin{aligned}
        \nu(M,v,g, \omega) & := \chi(W) - 3 \sigma(W) - \int_W c_1 c_3  ~ \in \mathbb{Z} / (48 \mathbb{Z}) \\
    \end{aligned}\end{equation}
\end{definition}
Here $\chi(W)$ and $\sigma(W)$ are the Euler characteristic and signature of $W$ respectively. 
Note the integral on the right hand side is sensible on $W$ since both $c_1$ and $c_3$ have compactly supported representatives.  
Proposition \ref{prop_u4_coboundary} proves the existence of such a coboundary, 
while Proposition \ref{prop_nu_u3} proves that the definition is well defined and independent of choice of coboundary. 

We define a further invariant in a similar fashion.
For manifolds with boundary the integral of non compactly supported classes is not well defined. 
However, integrating a product of classes is well defined
modulo some integer corresponding to the divisibility of its components. 

\begin{definition}\label{def_xi_u3}
    Let $M$ be as in Definition \ref{def_nu_u3}. 
    Additionally, let $m$ be the greatest divisor of $c_2 (M)$ modulo torsion
    and define the set 
    \begin{equation}\begin{aligned}
        S_m := \{ s \in H^4(M) \colon c_2(M) - ms \in TH^4(M) \} 
    \end{aligned}\end{equation}
    For $s \in S_m $ we define $u_s = c_2(W) - m\tilde{s} \in H^4(W; \mathbb{Q}) $ for some lift $\tilde{s} \in H^4(W) $ of $s$. 
    Let $\tilde{m} := \mathSet{lcm} ( m, 4) $.
    We define an invariant $\xi(M, g, \omega): S_m \rightarrow \mathbb{Q} / 3\tilde{m} \mathbb{Z} $ by
    \begin{equation}\begin{aligned}
        \xi(M, g, \omega)(s) & :=  7 \chi(W) - \frac{45\sigma(W) }{2} - \int_W \, \left( 7 c_{1} c_{3} - 2 c_{1}^{2} c_2 + \frac{c_{1}^{4} }{2} \right) 
        + \int_W \,\frac{3u_s ^{2}}{2} ~( \mathSet{mod}~ 3\tilde{m}) 
    \end{aligned}\end{equation}
\end{definition}
Proposition \ref{prop_xi_u3} proves that the definition is well defined and independent of choice of coboundary,
and of choices of lift $\tilde{s}$ for $s \in S$. 
In the case that $TH^4(M) = \{0 \}$, we write simply $\int_W c_2 ^2 $ for the ultimate term and consider $\xi \in (\mathrm{mod}~3\tilde{m})$ as a constant.

Although not the context which motivated this work, an obvious application of these invariants is to distinguish or classify almost contact $7$-manifolds.

\subsection{Application to holonomy $\lieGroup{G}_2$} 
In \cite{crowley2015new}, the authors define two invariants of 7-manifolds with $\lieGroup{G}_2$-structure that they denote by $\nu$ and $\xi$, 
and from which we have suggestively inherited our notation. 
Their invariants, which here we denote by $\nu'$ and $\xi'$, are defined via spin coboundary together with a spinor. 
According to \cite[Theorem 6.9]{crowley2015new}, for a $2$-connected closed $7$-manifolds $(\nu, \xi)$ determine the homotopy class of a $\lieGroup{G}_2$-structure of up to action of spin diffeomorphism.  

For a given 7-manifold with $\lieGroup{G}_2$-structure we say that a $\lieGroup{U}(3)$-structure with vanishing $c_3$ is compatible
if they share an $\lieGroup{SU}(3)$-reduction (Definition \ref{def_compatible_u3}). 
For a manifold with $\lieGroup{G}_2$-structure, a compatible $\lieGroup{U}(3)$-structure exists (Proposition \ref{prop_compatible_g2_u3_structures}). 

Proposition \ref{prop_su4_coboundary} proves the existence of a suitable $\lieGroup{SU}(4)$-coboundary to a 7-manifold with $\lieGroup{SU}(3)$-structure.
In doing so we are able to conclude the following result.

\begin{proposition}\label{prop_invariants_agree}
    Let $M$ be a closed spin $7$-manifold. 
    Suppose that $\varphi$ is a $\lieGroup{G}_2$-structure on $M$ and that $(g, \omega)$ is a compatible $\lieGroup{U}(3)$-structure on $M$. 
    Then $p_M = c_2(M)$ and 
    \begin{equation}\begin{aligned}
        \nu(M, g, \omega) & = \nu ' ( M, \varphi),  & \xi ( M, g, \omega) & = \xi' (M, \varphi)
    \end{aligned}\end{equation}
\end{proposition}

\subsection{The Twisted Connected Sum Construction} 

The TCS method has allowed for the construction of many $\lieGroup{G}_2$ manifolds \cite{kovalev2000twisted, corti2015}.
A very brief summary is given in Section \ref{sec_tcs_construction}. 
The construction starts from pairs of open manifolds with holonomy $\lieGroup{SU}(3)$ and tame asymptotic behaviour, 
called Asymptotically Cylindrical Calabi Yau manifolds (ACyl CYs) (Definition \ref{def_accy}). 
ACyl CYs are obtained from certain Kähler threefolds that have a K3 fibration, called building blocks (Definition \ref{def_building_block}).
Then $V := Z \setminus \Sigma$, for some K3 fibre $\Sigma \subset Z$, admits a Calabi-Yau structure (torsion free $\lieGroup{SU}(3)$-reduction) 
with non-compact end modelled on $\mathbb{R} \times S^1 \times \Sigma$ (see \cite[Theorem 3.4]{corti2015}). 
The Calabi-Yau structure asymptotes exponentially to a product cylinder, 
giving the K3 `at infinity' $\Sigma$ a hyper Kähler structure (ie torsion free $\lieGroup{SU}(2)$-reduction). 

The product $V \times S^1$ has the product $\lieGroup{G}_2$-structure. 
A pair $V_\pm \times S^1$ may be truncated and glued to obtain a $\lieGroup{G}_2$-manifold $(M, \varphi)$ 
(see \cite[Theorem 3.12]{corti2015}).
The gluing is constructed in part from a certain type of diffeomorphism called a \emph{Hyper Kähler Rotation} (Definition \ref{def_hkr}) $r:\Sigma_+ \rightarrow \Sigma_-$
between the K3s at infinity. 

For a TCS manifold we describe a compatible $\lieGroup{U}(3)$-structure to the $\lieGroup{G}_2$-structure in Section \ref{sec_tcs_coboundary}. 
We construct a suitable $\lieGroup{U}(4)$-coboundary which we call a TCS coboundary. 
The TCS coboundary allows for the computation of our invariants (See Proposition \ref{prop_tcs_nu_xi}). 

\begin{proposition}\label{prop_nu_is_24}
    For TCS manifold $(M, \varphi)$ the $\nu$-invariant $\nu(\varphi) = 24 ~(\mathSet{mod}~48)$.  
\end{proposition}
This agrees with \cite[Theorem 1.7]{crowley2015new}. 

\begin{proposition}\label{prop_xi_tcs_simplest} 
    Let TCS manifold $(M, \varphi)$ have TCS coboundary $W$. 
    Suppose that the cohomology of $M$ is torsion free, and that $m := \mathrm{gd}(p_M)$, the greatest divisor of the spin class $p_M \in H^4(M)$. 
    Then the $\xi$ invariant simplifies to $ \xi(\varphi) = \frac{3}{2} \int_W \, c_2 (W) ^{2}  ~( \mathSet{mod}~ 3\tilde{m}) $. 
\end{proposition}

\subsection{Examples} 

The TCS construction generates many examples of $\lieGroup{G}_2$-manifolds. 
For this class of $\lieGroup{G}_2$-manifolds, only the $\nu$-invariant had been computed previously 
but it is unable to distinguish homotopy classes of the resulting $\lieGroup{G}_2$-structure by Proposition \ref{prop_nu_is_24}.
We are now able to compute the $\xi$-invariant and use it to distinguish connected components of $\lieGroup{G}_2$-structures for TCS manifolds. 

The classification of 2-connected spin 7-manifolds \cite{crowley2014classification} up to spin diffeomorphism is given in terms of list of invariants. 
In particular, the spin diffeomorphism class of 2-connected 7-manifolds with torsion free cohomology, 
is determined by the triple of invariants $(b_3, m, \mu)$: the third Betti number; the greatest divisor of the spin class $m := \mathrm{gd}(p_M)$; and the generalized Eells-Kuiper invariant. 
In order to apply this, we restricted our choices in the TCS construction to those resulting in manifolds with torsion free cohomology and $b_2 (M)= 0$.  

We further note that for TCS manifolds, $m$ is an even divisor of $24$.
Hence that in the cases of torsion free cohomology,
the generalized Eells-Kuiper invariant reduces to a constant $\mu \in \mathbb{Z}/ \hat{m}\mathbb{Z}$, 
where $\hat{m} := \mathrm{gcd}\left(28, \mathrm{Num}\left(\frac{m}{4}\right) \right)$. 
Thus for TCS manifolds, $\mu(M)$ necessarily vanishes when $8 \nmid m $. 
In addition, the invariants $\mu, \nu$ and $\xi$ are subject to the relation (see \cite[Equation 38b]{crowley2015new})
\begin{equation}
    \label{eq_relating_mu_nu_xi}
    \frac{ \xi - 7 \nu }{12} = \mu ~\left(\mathrm{mod}~ \mathrm{gcd}\left(28, \frac{\tilde{m}}{4}\right)\right)
\end{equation}
from which we see that $\mu$ is completely determined by $\nu$ and $\xi$.
Similarly, if $6 \nmid m$ then $\xi$ is determined by $\mu$ and $\nu$. 
Thus we searched for pairs of TCS with torsion free cohomology and sharing $(b_3, m)$ with $m \in \{6, 12, 24\}$, and on finding them computed $\xi$. 

For a closed, $2$-connected spin $7$-manifold $M$ with torsion free cohomology, $M$ and $-M$ will have the same $(b_3, m)$ invariants. 
Moreover, in our context where $\hat{m} = 1$ or $2$, $\mu$ agrees on $M$ and $-M$. 
By the classification result then there exists an orientation reversing diffeomorphism $r: M \rightarrow M$.
A $\lieGroup{G}_2$-structure $\varphi$ is diffeomorphic to $-r^*(\varphi)$, yet $\xi(\varphi) = - (\xi(- r ^* \varphi))$. 
We look pairs of $\lieGroup{G}_2$-manifolds with torsion free cohomology, for which $(b_3, m)$ agree and $\xi$ differs and not simply by a sign change. 

The building blocks used in the TCS construction can be obtained from Fano threefolds.
Fano threefolds are well understood and much of the relevant data concerning them can be found in the literature. 
We considered examples of TCS manifolds involving the gluing of pairs of rank 1 and rank 2 Fano threefolds as was done in \cite{crowley2014exotic} to find exotic $\lieGroup{G}_2$ manifolds. 
However, this was insufficient for our needs---there were no 2-connected TCS manifolds with torsion free cohomology that shared $(b_3,m)$ and differing $\xi$-invariant.  
We continued by searching through matchings of rank 2 semi-Fanos.
Thanks to the work of \cite{arap2017existence} and others on which this builds, much of the data needed can also be read off available tables. 
\begin{example}
    The following two TCS manifolds have torsion free cohomology and are 2-connected with $(b_3, m)=(71,6)$. 
    As $\mu$ is vacuous so the underlying manifolds are diffeomorphic. 
    \begin{enumerate}
        \item The matching of two copies of the Fano threefold of Picard rank 1, index 1 and of genus $9$.
        \item A matching between a quadric $Q \subset \mathbb{P}^4$ blown up in smooth rational curve of degree 6, 
            and $\mathbb{P}^3$ blown up in a curve of degree 11 and genus 14. 
    \end{enumerate}
    The first has $\xi=0 ~ \mathrm{mod}~ 36$; while the second has $\xi=12 ~ \mathrm{mod}~ 36$. 
\end{example}

\begin{example}
    The following two TCS manifolds have torsion free cohomology and are 2-connected with $(b_3, m)=(85,24)$.
    \begin{enumerate}
        \item A matching of two copies of $\mathbb{P}^3$ blown up in a curve of degree 11 and genus 14. 
        \item A matching of $\mathbb{P}^1 \times \mathbb{P}^2$ and $\mathbb{P}^3$ blown up in a curve of degree 8 and genus 2. 
    \end{enumerate}
    The first has $\xi=12 ~ \mathrm{mod}~ 72$; while the second has $\xi=36 ~ \mathrm{mod}~ 72$. 
    By Equation \ref{eq_relating_mu_nu_xi}, $\mu = 1$ in both cases, so the underlying manifolds are diffeomorphic.
\end{example}

At a stage of the TCS construction one needs \emph{genericity results} (Definition \ref{def_genericity}). 
Roughly speaking, this guarantees that generic lattice polarized K3 appears as an anticanonical divisor in some building block, 
and ensures the existence of a gluing with the necessary properties. 
Although we have not solved systematically, we think a systematic approach is feasible. 
In doing would lead to many more examples. 
We have provided the genericity results relevant to our examples Section \ref{sec_examples}.

\subsection{Acknowledgements} 
I would like particularly to thank Johannes Nordstr\"om for many useful comments and suggestions. 
Thanks also to Jesus Martinez Garcia, and Diarmuid Crowley for assistance in understanding aspects of the algebraic geometry, 
and cobordism theory respectively
and to the EPSRC for their support.

\section{Background} 

\subsection{Certain representations of Lie groups} 
Representations form the linear algebra models of $G$-geometries. 
We recall some standard representations in terms of stabilizers of a general linear group in order to set out the notation and conventions used. 

Let $\{ e_i \} $ be the standard basis of $\mathbb{R} ^n$ and $\{e^i \}$ be its dual in $(\mathbb{R}^n )^*$.
The standard metric, volume form, symplectic form, complex structure, and complex volume form as follows. 
\begin{equation}\begin{aligned}
    g_{0,n} := & \textstyle \sum_{i=1} ^n e^i \otimes e^i  & 
    \mathrm{Vol}_{0,n} := & e^i \wedge \dots \wedge e^n  \\
    \omega_{0,m} := &\textstyle  \sum_{j=1} ^{m} e^{2j-1} \wedge e^{2j} &
    J_{0,m} := &\textstyle  \sum_{j=1} ^{m} (e^{2j-1} \otimes e_{2j} -  e^{2j} \otimes e_{2j-1})  \\
    \Omega_{0,m} := & \Lambda_{j=1} ^{m} (e^{2j-1} + i e^{2j})  \\
\end{aligned}\end{equation}
Then $\lieGroup{O}(n) = \mathrm{Stab}(g_{0,n})$, $\lieGroup{SO}(n) = \mathrm{Stab}(g_{0,n}, \mathrm{Vol}_{0,n})$. 
If $n=2m$ then $\lieGroup{U}(m) = \mathrm{Stab}(\omega_{0,m}, J_{0,m})$ and $\lieGroup{SU}(m) = \mathrm{Stab}(\omega_{0,m}, \Omega_{0,m})$.
We have 
\begin{equation}\begin{aligned}
    \lieGroup{SU}(m) \subset \lieGroup{U}(m) \subset \lieGroup{SO}(2m) \subset \lieGroup{O}(2m) 
\end{aligned}\end{equation}
and we note that $ \mathrm{Stab}(g_{0, 2m}, \omega_{0,m}) =  \mathrm{Stab}(g_{0, 2m}, J_{0,m}) =  \mathrm{Stab}(J_{0,m}, \omega_{0,m})$. 
For $n = 2m+k$, it is convenient to view $\lieGroup{U}(m) \cong \mathrm{Stab}(g_{0,n}, \omega_{0,m}, e_{m+1}, \dots, e_{m+k})$, 
and likewise for $\lieGroup{SU}(m)$. 
When the dimension of the ambient space is clear we shall omit it from the notation. 

We have two additional forms
\begin{equation}\begin{aligned}
    \varphi_0 := & e^{123} + e^{145} + e^{167} + e^{246} - e^{257} - e^{347} - e^{356} \in \Lambda^3 ( \mathbb{R} ^7 )^*  \\
   \psi_{0} := & e^{1234} + e^{1256} + e^{1278} + e^{1357} - e^{1368}
    -e^{1458} - e^{1467} - \\ &  e^{2358} - e^{2367} - e^{2457} + e^{2468}
    + e^{3456} + e^{3478} + e^{5678} \; \in \Lambda^{4}(\mathbb{R}^{8})^{*}
\end{aligned}\end{equation}
where $e^{ij} = e^i \wedge e^j$ etc. 
The first form is the standard $\lieGroup{G}_2$-form since $\lieGroup{G}_2 = \mathrm{Stab}(\varphi_0)$. 
We find that $\lieGroup{G}_2 \subset \lieGroup{SO}(7) $ (see Bryant \cite[section 2]{bryant1987metrics}) and it acts transitively on the sphere $S^6$. 
The stabilizer of any unit vector is isomorphic to $\lieGroup{SU}(3)$.
The second form corresponds to the image of the spinor representation $\Delta : \lieGroup{Spin}(7) \rightarrow \lieGroup{SO}(8)$.
That is $\Delta(\lieGroup{Spin}(7)) = \mathrm{Stab}(\psi_0)$. 
$\mathrm{Stab}(\psi_0)$ acts transitively on the sphere $S^7$, with the stabilizer of any unit vector isomorphic to $ \lieGroup{G}_2$. 

For groups $\lieGroup{Spin}(n)$ and $\lieGroup{Spin}^c(n)$ we have a slightly different approach.
Recall that $\lieGroup{Spin}(n) \rightarrow \lieGroup{SO}(n)$ 
is a double cover, 
and that 
$\lieGroup{Spin}^c(n) := \lieGroup{Spin}(n) \times _{\mathbb{Z}_2} \lieGroup{U}(1)$ 
is a double cover $\lieGroup{Spin}^c (n) \rightarrow \lieGroup{SO}(n) \times \lieGroup{U}(1)$. 

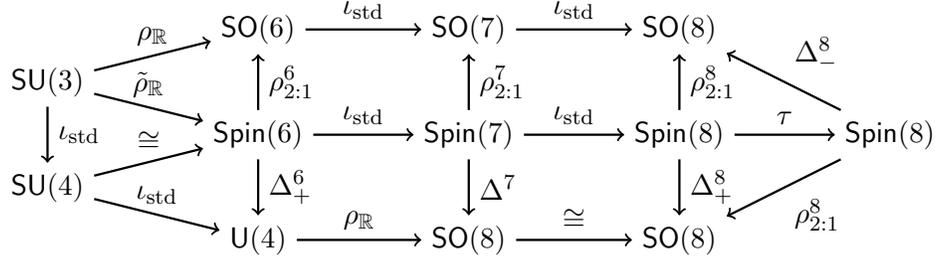
\begin{figure}
    \begin{tikzpicture}[scale=0.7]
        \tikzstyle{label} = []
        \node[label] (00) at (0,0) {$ \lieGroup{SO}(6) $}; 
        \node[label] (01) at (4,0) {$ \lieGroup{SO}(7) $}; 
        \node[label] (02) at (8,0) {$ \lieGroup{SO}(8) $}; 
        \node[label] (10) at (0,-2) {$ \lieGroup{Spin}(6) $}; 
        \node[label] (11) at (4,-2) {$ \lieGroup{Spin}(7) $}; 
        \node[label] (12) at (8,-2) {$ \lieGroup{Spin}(8) $}; 
        \node[label] (13) at (12,-2) {$ \lieGroup{Spin}(8) $}; 
        \node[label] (1m1) at (-4,-1) {$ \lieGroup{SU}(3) $}; 
        \node[label] (2m1) at (-4,-3) {$ \lieGroup{SU}(4) $}; 
        \node[label] (20) at (0,-4) {$ \lieGroup{U}(4) $}; 
        \node[label] (21) at (4,-4) {$ \lieGroup{SO}(8) $}; 
        \node[label] (22) at (8,-4) {$ \lieGroup{SO}(8) $}; 
        \draw[->, line width = 0.3mm] (00) -- (01) node[midway, above] {$ \iota_{\mathrm{std}} $};
        \draw[->, line width = 0.3mm] (01) -- (02) node[midway, above] {$ \iota_{\mathrm{std}} $};
        \draw[->, line width = 0.3mm] (10) -- (11) node[midway, above] {$ \iota_{\mathrm{std}} $};
        \draw[->, line width = 0.3mm] (11) -- (12) node[midway, above] {$ \iota_{\mathrm{std}} $};
        \draw[->, line width = 0.3mm] (12) -- (13) node[midway, above] {$ \tau $};
        \draw[->, line width = 0.3mm] (1m1) -- (00) node[midway, above] {$ \rho_\mathbb{R} $};
        \draw[->, line width = 0.3mm] (1m1) -- (10) node[midway, above] {$ \tilde{\rho}_\mathbb{R} $};
        \draw[->, line width = 0.3mm] (2m1) -- (10) node[midway, above]  {$ \cong $};
        \draw[->, line width = 0.3mm] (2m1) -- (20) node[midway, above] {$ \iota_{\mathrm{std}} $};
        \draw[->, line width = 0.3mm] (20) -- (21) node[midway, above] {$ \rho_\mathbb{R} $};
        \draw[->, line width = 0.3mm] (21) -- (22) node[midway, above] {$ \cong $};
        \draw[->, line width = 0.3mm] (1m1) -- (2m1) node[midway, right] {$ \iota_{\mathrm{std}} $};
        \draw[->, line width = 0.3mm] (10) -- (00) node[midway, right] {$ \rho_{2:1} ^6 $};
        \draw[->, line width = 0.3mm] (10) -- (20) node[midway, right] {$ \Delta_+ ^6 $};
        \draw[->, line width = 0.3mm] (11) -- (01) node[midway, right] {$ \rho_{2:1} ^7 $};
        \draw[->, line width = 0.3mm] (11) -- (21) node[midway, right] {$ \Delta ^7 $};
        \draw[->, line width = 0.3mm] (12) -- (02) node[midway, right] {$ \rho_{2:1} ^8 $};
        \draw[->, line width = 0.3mm] (12) -- (22) node[midway, right] {$ \Delta_+ ^8 $};
        \draw[->, line width = 0.3mm] (13) -- (02) node[midway, above right] {$ \Delta_- ^8 $};
        \draw[->, line width = 0.3mm] (13) -- (22) node[midway, below right] {$ \rho_{2:1} ^8 $};
    \end{tikzpicture}
    \caption{A commutative diagram of relevant Lie group homomorphisms} 
    \label{fig_lie_group_hom} 
\end{figure}

It is not critical although may be of assistance to some readers to note the following relations between representations. 
The half spin representation $\Delta_+ ^6: \lieGroup{Spin}(6) \rightarrow \lieGroup{U}(4)$ is faithful with image $\lieGroup{SU}(4)$. 
There is a `triality' property of $\lieGroup{Spin}(8)$ that permutes the three $\lieGroup{SO}(8)$-representations: the double cover $\rho_{2:1} ^8$; 
and the two half spin representations $\Delta _\pm ^8$. See Figure \ref{fig_lie_group_hom}.

\subsection{$G$-structures} 
\label{G_structures} 

Let $G \subset H$ be a closed Lie subgroup. 
Then the space of left cosets $F = H/G$ is naturally a homogeneous $H$-space. 
If $X$ is a manifold and $E \rightarrow X$ a principal $H$-bundle on $X$, then $F \times_H E$ is a smooth $F$-fibre bundle over $X$. 
In greater generality, suppose Lie group $H$ acts on a space $K$, fix some $k \in K$ and let $G := \mathSet{Stab}_H (k)$ be the stabilizing subgroup. 
Let us call a section $s_k \in \Gamma (K \times_H E)$ \emph{$k$-like} if for each $x \in X$ there exists some $e \in E_x $ such that $(e,k) \in s_k (x) $.
We can ask whether there exists a $k$-like section of $K \times_H E$.
By taking $F := \mathrm{Orb}_H (k)$, the orbit of $k \in K$ by $H$, and considering the action of $H$ on $F$ we have a transitive action. 
Provided that the stabilizer $G := \mathrm{stab}_{H}(k)$ of $k$ is closed, then $F \cong H/G$. 
If the action of $H$ is transitive then any section is $k$-like. 

\begin{lemma}
    Suppose that $H$ is a Lie group, $E \rightarrow X$ is a principal $H$-bundle on a manifold $X$, and space $K$ admits an $H$-action. 
    Suppose we fix some $k \in K$ and that $G := \mathSet{Stab} (k) \subset H$ is closed. 
    There is a one-to-one correspondence between $k$-like sections $K \times_H E$ and $G$-structures on $E$.  
\end{lemma}

In light of this, a section that uniquely determines a $G$-structure may itself be referred to as a $G$-structure. 
For example:
\begin{enumerate}
    \item A $\lieGroup{G}_2$-structure on a $X^7$ is a 3-form $\varphi \in \Gamma(\Lambda^3 T^* X)$ that is $\varphi_0$-like. 
    \item A $\lieGroup{U}(m)$-structure on a $X^{2m}$ is a section $(\omega, J) \in \Gamma(\Lambda^2 T^* X \otimes \mathrm{End}(TX))$ that is $(\omega_0, J_0)$-like. 
    \item An $\lieGroup{SU}(3)$-structure on a $X^{7}$ is a section $(v, \omega, \Omega) \in \Gamma(TX \otimes \Lambda^2 T^* X \otimes \Lambda^m _\mathbb{C} (T^*X))$ that is $(e_1,\omega_0, \Omega_0)$-like.
    \item A $\lieGroup{Spin}(7)$-structure on a $X^{8}$ is a section $\psi \in \Gamma(\Lambda^4 T^* X)$ that is $\psi_0$-like.
\end{enumerate}
and so on. 

In case of spin and spin$^c$, the story is slightly different. 
Let a manifold $X^n$ have an $\lieGroup{SO}(n)$-structure denoted by $F_{\lieGroup{SO}} X $. 
A spin structure $E \rightarrow X$ is a principal $\lieGroup{Spin}(n)$-bundle together with a double cover $E \rightarrow F_{\lieGroup{SO}} X$, 
such that the standard double cover representation $\rho: \lieGroup{Spin}(n) \rightarrow \lieGroup{SO}(n)$, commutes with respective group actions.
A spin manifold will mean a manifold with a given spin structure.  
Let $X$ be equipped with a Hermitian line bundle $L$. 
A spin$^c$ structure $E \rightarrow X$ on $(X,L)$ is a principal $\lieGroup{Spin}^c(n)$-bundle together with a double cover $E \rightarrow F_{\lieGroup{SO}}X \otimes F_{\lieGroup{U}}L$, 
such that the standard double cover representation $\rho: \lieGroup{Spin}^c (n) \rightarrow \lieGroup{SO}(n) \times \lieGroup{U}(1)$ commutes with respective group actions.
We may refer to $L$ as the associated line bundle of the spin$^{c}$ structure.
We may refer to $E$ as spin$^c$ structure on $X$ if the associated $L$ is implicit or unimportant. 

In the case that a manifold $X$ is \emph{a priori} endowed with an $H$-structure, and $G \subset H$, then a $G$-structure on $X$ will be assumed to factor through the $H$-structure. 
For example, if $X^7$ is a spin manifold, any $\lieGroup{G}_2$-structure on $X$ will induce the metric and orientation agreeing with that of the spin structure on $X$;
an $\lieGroup{SU}(3)$-structure on a manifold $X^7$ with $\lieGroup{G}_2$-structure is defined by a nowhere vanishing vector field and so on. 

The embedding $\lieGroup{G}_2 \hookrightarrow \lieGroup{Spin}(7)$ induces a map from a $\lieGroup{G}_2$-structure to a spin structure. 
A $\lieGroup{U}(m)$-structure induces a canonical spin$^{c}$ structure. 
Explicitly, the map $ i \oplus \mathSet{det} : \lieGroup{U}(m) \rightarrow \lieGroup{SO}(2m) \times \lieGroup{U}(1) $
where $i$ is the standard embedding forgetting complex structure, lifts to spin$^{c}$. 
The associated line bundle is then the $\mathSet{det}$ bundle associated to the $\lieGroup{U}(m)$-structure.
More generally, this maybe applied to $n$-manifolds provided that $2m\leq n$.  
We now make precise what we meant by a manifold having compatible $\lieGroup{G}_2$-structure and $\lieGroup{U}(3)$-structure. 

\begin{definition}\label{def_compatible_u3}
    Let $M$ be a spin 7-manifold with both $\lieGroup{G}_2$-structure $\varphi$. 
    A $\lieGroup{U}(3)$-structure $(\omega, J)$ on $M$ such that $c_3 = 0 $ is said to be \emph{compatible} with $\varphi$
    if there exists $\lieGroup{SU}(3)$-reductions of both $(\omega, J)$ and $\varphi$ which are identical 
    when considered as reductions of the $\lieGroup{SO}(7)$-structure. 
\end{definition}

\begin{definition}\label{def_structured_coboundary} 
    Let $M$ be a closed $n$-manifold with $G$-structure $E \rightarrow F_{\lieGroup{SO}}M$. 
    Suppose that $H \subset \lieGroup{SO}(n+1)$ acts transitively on the $n$-sphere and $\lieGroup{Stab}_H(e_{n+1}) \cong G$. 
    Then an $H$-coboundary $W$ to $M$, is a coboundary $W$ together with an $H$-structure $E_H \rightarrow F_{\lieGroup{SO}} W$
    such that $\partial W \cong M$ induces $E \rightarrow E_H$ as a $G$-structure. 
\end{definition}
This definition can be sensibly extended for spin and spin$^c$, 
as well as further generalized for stable structured coboundaries.

\begin{definition}\label{def_homotopic_structures}
    Let $H$ be a Lie group, $X$ a manifold and $E \rightarrow X$ a principal $H$-bundle on $X$. 
    Let $G \subset H$ be a Lie subgroup, and let $F = H/G$. 
    Two $G$-structures $s_i \in \Gamma(E_F)$ are said to be homotopic if there exists homotopy
    \begin{equation}
        S : I \times X \rightarrow E_F 
    \end{equation}
    such that $S(0, x) = s_0 (x)$ and $S(1, x) = s_1 (x)$, and for all $t \in I$, $S(t, \cdot) \in \Gamma(E_F)$. 
\end{definition}

Note that if $K \subset G \subset H$ are Lie groups, two $K$-structures of principal $H$-bundle $E \rightarrow X$ maybe homotopic, 
but not when considered as reductions of $G$-structures $P \rightarrow X$.

\subsection{Structure reductions} 

In proving the existence of coboundaries with certain properties of the structure group, 
we require `improving' an $H$-structure on some manifold $X$ to a $G$-reduction. 
We have seen that this is equivalent to finding a section of the associated $F$-fibre bundle, where $F=H/G$. 
We employ a standard obstruction theory technique inducting extensions on skeleta. 

\begin{proposition}\cite[Section 3.3]{hatcher2018vector}
    Let $(X,A)$ be a finite CW pair. 
    Suppose that $E \rightarrow X$ is an $F$-fibre bundle over $X$. 
    We will assume that the base space $X$ is connected; 
    that $F$ is path-connected and $\pi_1 (F)$ acts trivially on $\pi_n (F)$ for all $n$;
    and that $\pi_1 (X)$ acts trivially on $\pi_n (F)$. 

    Suppose that $s \in \Gamma(E|_A)$ is a section of $E$ over $A$. 
    If $F$ is $(k-1)$-connected, then the first non-trivial obstruction space is $H^{k+1} (X, A; \pi_k (F))$.
    Moreover, the obstruction class $c(E,s) \in H^{k+1} (X,A ; \pi_k (F))$ of $s$ is independent of choices.
\end{proposition}
In other words, the primary obstruction $c(E, s) \in H^{k+1}(X,A; \pi_k(F))$ is natural (ie functorial). 

In applications such as our own the fibre bundle is derived from a quotient of a principal bundle. 
Thus the fibre is a homogeneous space. 
Suppose that $ \sigma: H \rightarrow \lieGroup{SO}(k+1)$ is a real representation of a compact Lie group $H$ for $k$ odd.
The associated vector bundle $E_V := \mathbb{R}^{k+1} \times_\sigma E$ is oriented. 
Assume that $H$ acts transitively on the sphere $S^{k} \subset V$. 
Let $G := \mathrm{Stab}_H (v) $, the stabilizer of a unit vector $ v \in S^k$. 
Suppose $s_A \in \Gamma(E_S|_A)$ is a section of the associated sphere bundle $E_S = S^k \times _\sigma E$ over $A$.
Define $A' := \mathrm{Image}(s_A) \subset E_S$. 
We have a projection $p:(E_V,A') \rightarrow (X,A)$, and embedding $j: (E_V,A') \rightarrow (E_V,E_S)$. 
Note that $p$ is a retract and hence induces an isomorphism at the level of cohomology. 
Let $\tau \in H^{k+1}(E_V, E_S)$ be the Thom class. 
The Euler class relative to $s$ is defined as $e(E_V,s) := ((p^*)^{-1}\circ j^*) (\tau)$. 
In the case that $A$ is empty, then we have the absolute Euler class in the commonly understood sense and denoted $e(E_V)$. 

Sharafutdinov \cite[Theorem 1.1]{sharafutdinov1973relative} claims this is the unique class satisfying the following axioms.
\begin{enumerate}
    \item Naturality: $f^* e(E_V,s) = e ( f^* E_V, f^* s)$ for a morphism of CW-pairs $f:(Y,B) \rightarrow (X,A)$. 
    \item Multiplicative: $ e(E_V \oplus E'_{V'}, s) = e(E_V, s) e(E' _{V'})$ for even dimensional oriented vector bundles $E_V$, and $E'_{V'}$, and section $s \in \Gamma(E_S |_A)$. 
    \item Normed: If $E_V = \mathcal{O}_{\mathbb{P}^1}(1)$ as a real vector bundle, and $X = S^2 = \mathbb{P}^1$, 
        then $e(E_V) \in H^2(X)$ is the oriented generator. 
\end{enumerate}

\begin{proposition}\label{prop_euler_class_obstruction} 
    With the notation and assumptions as above, the relative Euler class is the primary obstruction to extending $s_A \in \Gamma(E_S |_A)$ to $X$. 
    That is 
    \begin{equation}
        c(E_S, s) = e(E_V, s) 
    \end{equation}
\end{proposition}

For connected manifolds $X$ and $X'$ (either with or without boundary), 
if $X$ and $X'$ are oriented, then their connected sum $X \# X'$ is oriented. 
Moreover, if $X$ and $X'$ are spin manifolds, then $X \# X'$ is spin and is spin cobordant to $X \sqcup X'$
(see \cite[Lemma 2.1]{kervaire1963groups}). 
Thus we have the following formula for spin 8-manifolds.
\begin{proposition}
    \label{prop_connected_sum_euler_classes} 
    Let $X,X'$ be closed spin 8-manifolds. 
    \begin{equation}\begin{aligned}
        e_+ ( X\# X') = e_+ (X) + e_+ (X') -1  
    \end{aligned}\end{equation}
    Furthermore, if $X,X'$ be are spin 8-manifolds with boundary
    \begin{equation}\begin{aligned}
        e_+ ( X\# X', s \sqcup s') = e_+ (X, s) + e_+ (X', s') -1  
    \end{aligned}\end{equation}
    where $s,s'$ are unit spinors of the boundaries of $X$ and $X'$ respectively.
\end{proposition}
Thus by taking connected sums, we can effectively kill off the relative Euler class to ensure the existence of, \emph{eg} a nowhere vanishing spinor.

\section{The Invariants} 

We first prove the existence of certain coboundaries. 
We use this to prove the well definedness of $\nu$ and $\xi$.

\subsection{Existence of certain coboundaries} 

\begin{proposition}\label{prop_su4_coboundary}
    For any closed manifold $M^7$ with $\lieGroup{SU}(3)$-structure 
    there exists a $\lieGroup{SU}(4)$ coboundary $W$ to $M$ over $H^4$. 
\end{proposition}

\begin{proof}
    An $\lieGroup{SU}(3)$-structure on a 7-manifold determines a spin structure and a pair of orthonormal unit spinors. 
    That is, the $\lieGroup{SU}(3)$-structure is precisely the reduction of the spin structure defined as the stabilizer of the orthonormal pair of spinors. 

    According to \cite[pg 204]{milnor1963spin}, the spin cobordism group $\Omega_{\lieGroup{Spin}} ^7 = 0$.
    Thus a spin coboundary of spin manifold $M$ exists. 
    Moreover, \cite[Theorem 3]{milnor1961killing} we can choose that the coboundary $W$ is 3-connected.  
    By Poincare-Lefschetz duality $H^5(W,M) \cong H^3(W) = 0$. 
    The long exact sequence of cohomology given by the pair $(W,M)$, then implies $H^4(W) \rightarrow H^4(M)$ is onto. 
    That is $W$ is a coboundary over $H^4$. 

    Let $s_i \in \Gamma (S_M)$ be the unit spinors determined by the $\lieGroup{SU}(3)$-structure. 
    It remains to extend the unit spinors $s_i \in \Gamma(S_M)$ to the interior of $W$.
    By Proposition \ref{prop_euler_class_obstruction}, the primary obstruction to extending $s_1$ to skeleta of $W$ is the relative Euler class $e(S_W, s_1) \in H^8(W, M; \pi_7(S^7)) \cong \mathbb{Z}$. 
    Note that $e_+(S^4 \times S^4) = 2$, and $e_+(\mathbb{HP}^2) = 0$. 
    By Proposition \ref{prop_connected_sum_euler_classes}, taking successive connected sums of $W$ with one of these, 
    we may assume that $e_+(W, s_1) = 0 $ and so admits a nowhere vanishing spinor field extending $s_1$. 
    Note that $W$ remains three connected.

    Now we consider extending $s_2$ such that it remains perpendicular to the extension of $s_1$. 
    This is equivalent to extending the section of an $S^6$-fibre bundle. 
    $W$ is simply connected, so $H^7(W,M; \pi_6 (S^6)) = 0$. 
    Thus the primary obstruction of extending $s_2$ vanishes, and we have an extension of $s_2$ to $W^7$. 
    The secondary obstruction space is $H^8(W,M; \pi_7 ( S^6) ) \cong \mathbb{Z}_2$. 
    We consider extending $s_2$ at the level of cells. 
    There are precisely two possible homotopy classes: 
    $s_2$ over $\partial B_j ^8 $ corresponds to $ 0 \in \pi_7 (S^6)$, 
    and $s_2$ corresponds to $1 \in \pi_7 (S^6)$. 
    In the former case, we can extend $s_2$ over $B_j ^8$, while in the latter we can not. 
    
    If for each $8$-cell we can extend $s_2$ from the boundary $S^7$ to the interior, then we are done. 
    Assume there exists a cell $(B^8, S^7) \rightarrow (W^8, W^7)$ for which this is not possible. 
    We construct now an $8$-manifold with which we can replace each offending 8-cell and over the resulting manifold can extend $s_2$. 

    The manifold $\mathbb{HP}^2 $ is spin and admits a $\lieGroup{Spin}(7)$-structure but not an $\lieGroup{SU}(4)$-structure \cite[Theorem 5.7]{bowden2014stein}. 
    Moreover, puncturing at some $p\in \mathbb{HP}$ to obtain $\mathbb{HP}^2 \setminus \{p\} $ the $\lieGroup{Spin}(7)$-structure does admit an $\lieGroup{SU}(4)$-reduction.
    Fix an $\lieGroup{SU}(4)$-structure $E \rightarrow \mathbb{HP}^2 \setminus \{ p \}$. 
    The restriction of $E$ to an $S^7$ boundary of a neighbourhood of $p \in \mathbb{HP}$ is equivalent to that of a troublesome 8-cell. 
    
    The Euler characteristic $\chi(\mathbb{HP}^2) = 3$, while introducing three punctures results in $\chi (\mathbb{HP}^2 \setminus \{p, p', p'' \}) = 0$. 
    Thus there exists a nowhere vanishing vector field assumed normal on the boundary. 
    Fix a nowhere vanishing vector field and define an orientation reversing bundle involution $i: T(\mathbb{HP}^2 \setminus \{p, p', p''\}) \rightarrow  T(\mathbb{HP}^2 \setminus \{p, p', p''\})$. 
    We restrict $E$ to $\mathbb{HP}^2 \setminus \{ p, p', p'' \} $. 
    In a $B^8$ neighbourhood of $p'$ and $p''$ we fix a trivialization of $E$.
    Let $E'$ be the $\lieGroup{SU}(4)$-structure induced via $i$.
    $T^8$ has the flat $\lieGroup{SU}(4)$-structure, and we fix a $B^8$ subset. 
    Thus we can glue together an $\lieGroup{SU}(4)$-structure on $(-\mathbb{HP}^2)\# T^8 \# T^8 \setminus \{p \}$ agreeing with $E'$ on restriction.
    Let $X = (-\mathbb{HP}^2) \# T^8 \# T^8$. 
    Then $X \setminus \{p \}$ is an $8$-manifold with boundary $S^7$ that has an $\lieGroup{SU}(4)$-structure
    which on restriction to the boundary agrees with that of the offending $8$-cells. 
    Thus, we replace each such $8$-cell with $X \setminus \{ p \}$. 

    Let $W'$ be the resulting manifold after this process. 
    We retain the property that $H^4(W') \rightarrow H^4(M) $ is onto. 
\end{proof}

\begin{proposition}\label{prop_u4_coboundary} 
    For any closed manifold $M^7$ with $\lieGroup{U}(3)$-structure 
    there exists a $\lieGroup{U}(4)$ coboundary $W$ to $M$ over $H^4$. 
\end{proposition}

\begin{proof}
    The $\lieGroup{U}(3)$-structure on $M$ defines a spin$^c$ structure on $M$. 
    The stable unitary cobordism group $\Omega_\lieGroup{U} ^7 = 0$ (see, for example, Stong \cite[Chapter IV, Example 3]{stong1973cobordism} which summarises the results of Milnor \cite{milnor1960cobordism} and Novikov \cite{novikov1960problems, novikov1962homotopy}). 
    Thus there exists a stable unitary coboundary $W$ to $M$.  
    The unitary structure on $W$ defines a $\lieGroup{Spin}^c(8+N)$-structure on $TW \oplus \mathbb{R}^N$. 
    Restricting to frames of $TW$, we recover a spin$^c$ structure on $TW$.
    Moreover, $W$ is a spin$^{c}$-coboundary to spin$^c$ manifold $M$. 
    By the same arguments of surgery theory used in the spin coboundary case, we can assume that the homotopy $W$ agrees with that of $B\lieGroup{Spin}^c$ below the middle dimension.
    In particular we can assume that $H^3(W) = 0$, so that $\pi(W) = 0$ and $H^4(W) \rightarrow H^4(M)$ is onto. 

    We improve the spin$^{c}$ structure on $W$ to a $\lieGroup{U}(4)$-structure agreeing with the $\lieGroup{U}(3)$-structure on the boundary. 
    As 
    \begin{equation}
        \lieGroup{Spin}^{c}(8) / \lieGroup{U}(4) \cong \lieGroup{Spin}(8) / \lieGroup{SU}(4) \cong S^7 \times S^6 
    \end{equation}
    the proof proceeds identically to Proposition \ref{prop_su4_coboundary}.  
\end{proof}

Note that while going via a spin$^c$-coboundary may seem a like a detour, 
it is not obvious how one would instead improve a stable unitary structure to a genuine $\lieGroup{U}(4)$-structure in an analogous manner.

\subsection{The $\nu$-invariant} 

The $\nu$ invariant is the boundary defect of a combination of two characteristic class formulas for closed almost complex 8 manifolds: the Hirzebruch-Riemann-Roch theorem, and the Hirzebruch signature theorem. 
We recall the index theorem for the spin$^{c}$ Dirac operator \cite[Theorem D.15]{lawson1989spin}
which we wish to express in terms of Chern classes. 

\begin{lemma}
    Let $X$ be a closed $2n$-manifold with spin$^{c}$ structure and let $L$ be the associated $\mathbb{C}$-line bundle.
    Let $D: \Gamma S^+(X) \rightarrow \Gamma S^- (X)$ be the spin$^{c}$ Dirac operator. 
    Then 
    \begin{equation}
        \mathSet{ind} D = \int_X \hat{A}(X) e^{z /2} 
    \end{equation}
    where $z = c_1(L)$.
\end{lemma}

The expansion of $\hat{A}$ up to terms of degree at most $8$ is 
\begin{equation}
    1 -\frac{5 p_1 }{2^3 \cdot 3} + \frac{- 4 p_2 + 7 p_1 ^2 }{2^7 \cdot 45} 
\end{equation}
while $e^{z/2}$ expands to
\begin{equation}
    1 + \frac{z}{2} + \frac{z^{2}}{2^3} + \frac{z^{3}}{2^4 \cdot 3} + \frac{z^{4}}{2^7 \cdot 3} 
\end{equation}
In the case that $X$ has a $\lieGroup{U}(n)$-structure then $p_1 = -2c_2 + c_1 ^2 $, and $p_2 = 2c_4 - 2c_1 c_3 + c_2 ^2$, 
while the associated line bundle is the determinant bundle of the almost complex tangent space, 
so $z = c_1(X)$. 
Thus, the index theorem in terms of Chern classes is
\begin{equation}\label{eq_HRR} 
    720 \cdot \mathSet{ind}(D) = \int_X -c_{1}^{4} + 4 \, c_{1}^{2} c_{2} + 3 \, c_{2}^{2} + c_{1} c_{3} - c_{4}
\end{equation}
For a genuine complex manifold $\mathrm{ind}(D)$ is the holomorphic Euler characteristic, 
and equation (\ref{eq_HRR}) is the Hirzebruch-Riemann-Roch theorem in the almost complex setting (see Hirzebruchs ICM address on complex manifolds \cite[pg 122]{hirzebruch1958komplexe}).

The Hirzebruch signature theorem for closed 8 manifolds states that 
\begin{equation}
    \sigma(X) = \frac{1}{45} \int_X ( 7 p_2 - p_1 ^2) 
\end{equation}
In terms of Chern classes this becomes 
\begin{equation}\label{eq_hizebruch_chern}
    45 \cdot \sigma(X) = \int_X -c_{1}^{4} + 4 \, c_{1}^{2} c_{2} + 3 \, c_{2}^{2} - 14 \, c_{1} c_{3} + 14 \, c_{4}
\end{equation}

We make the following note on cohomology. 
There are several cup products available in the case of relative cohomology (see \cite[section 3.2]{hatcher2002algebraic}).
For a ring $R$ and topological pair $(A,B)$, let $ H^* _0 (A;R) := \mathSet{Im}(H^*(A,B;R) \rightarrow H^* (A;R)) $. 
There exists a natural graded product
\begin{equation}\begin{aligned}
    \label{eq_cup_prod_on_H0}
    H^i _0 (A;R) \times H^j _0 (A;R) & \rightarrow H^{i+j} (A, B;R)
\end{aligned}\end{equation}
which can be defined by lifting one of the factors to $H^*(A,B;R)$, and considering the product structure on $H^i(A,B;R) \times H^j(A;R)$. 
In particular, for a compact $8$-manifold with boundary $(W,M)$ we can sensibly integrate degree 8 classes of $H^*(W;R)$ that are given as products of classes in the image of $H^*(W,M;R)$. 

\begin{proposition}\label{prop_nu_u3}
    Definition \ref{def_nu_u3} of the $\nu$ invariant is well defined and independent of choice of coboundary $W$. 
\end{proposition}

\begin{proof} 
    Comparing eqs (\ref{eq_HRR}) and (\ref{eq_hizebruch_chern}), we eliminate the first three terms of the integrand.
    Integrating the top Chern class $c_4(X) $ agrees with the Euler characteristic $\chi(X)$. 
    Thus 
    \begin{equation}
        48 \cdot \mathSet{ind}(D) = -\chi(X) + 3 \sigma(X) + \int_X c_{1} (X)  c_{3} (X) 
    \end{equation}
    In particular, the right hand side is a multiple of $48$.

    Suppose that $W_0$ and $W_1 $ are $\lieGroup{U}(4)$-coboundaries of $M$ as in Proposition \ref{prop_u4_coboundary}. 
    Let $-M$ be the manifold identical to $M$, but with $\lieGroup{U}(3)$-structure $(-v, g, \omega)$ and hence reversed orientation. 
    Let $W_2$ be a $\lieGroup{U}(4)$ coboundary to $-M$. 
    We may form two closed 8 manifolds $X_i = W_i \cup_{M} W_2 $, for $i=0,1$, and the $\lieGroup{U}(4) $-structures of their components can be glued to give $\lieGroup{U}(4)$-structures to $X_i$. 

    The Euler characteristics split $ \chi(X_i) = \chi(W_i) + \chi(W_2)$. 
    By Novikov additivity \cite[Proposition 7.1]{atiyah1968index}, the signatures also split as $\sigma(X_i) = \sigma(W_i) + \sigma(W_2) $. 
    Since $c_1(M)$ and $c_3(M)$ vanish, so $\int_X c_1 (X)c_3(X) = \int_{W_i} c_1(W_i) c_3(W_i) + \int_{W_2} c_1(W_2) c_3(W_2) $.
    Thus
    \begin{equation}
         \chi(W_0) - 3 \sigma(W_0) - \int_{W_0} c_{1} c_{3} 
        =  \chi(W_1) - 3 \sigma(W_1) - \int_{W_1} c_{1} c_{3} ~ (\mathSet{mod} ~ 48) 
    \end{equation}
    It follows that $\nu(M)$ is independent of choice of coboundary. 
\end{proof}

\subsection{The $\xi$-invariant} 

For a complex vector bundle $E \rightarrow M$ over a closed manifold $M$, $c_{k}(E) = w_{2k}(E) ~\mathrm{mod}~2$, 
where $w_k(E)$ is the $k^\mathrm{th}$ Stieffel Whitney class of $E$. 
The $4^\mathrm{th}$ Wu class $v_4(M) = w_4(M) + w_2 ^2(M)$, which necessarily vanishes on a manifold of dimension $\leq 7$.
Thus for a 7-dimensional manifold $M$ with $\lieGroup{U}(3)$-structure, $c_2(M) + c_1(M) ^2 = 0 ~\mathrm{mod}~2$.
In particular, if $c_1(M) ^2 = 0$, then $2|c_2(M)$.  
It follows that $\mathSet{lcm}(m,4)$ is either $m$ or $2m$.

\begin{proposition}\label{prop_xi_u3}
    Definition \ref{def_xi_u3} of the $\xi$ invariant is well defined and independent of choice of coboundary $W$. 
\end{proposition}

\begin{proof}
    We argue analogously to the above, but starting from the Hirzebruch signature theorem. 
    Let $X$ be a closed $8$ manifold with $\lieGroup{U}(4)$-structure. 
    Then eq. (\ref{eq_hizebruch_chern}) says that
    \begin{equation}
        0 = 14 \chi(X) - 45 \sigma(X) +  \int_X -c_{1}^{4} + 4 \, c_{1}^{2} c_{2} + 3 \, c_{2}^{2} - 14 \, c_{1} c_{3}
    \end{equation}

    Let $W_i$'s be as above, and $X_i = W_i \cup_M W_2$. 
    We have the splitting results for $\chi$ and $\sigma$, 
    and all terms in the integrand are well defined on $W_i$, bar $c_2 ^2$, via Equation \ref{eq_cup_prod_on_H0}.  
    We require a splitting result in order to make sense of $c_2(X) ^2$. 

    For $x \in H^4(X; \mathbb{Q})$, $\int_X c_2(X)^2 = \int_X (c_2 (X) - mx)^2 ~\mathrm{mod} ~2\tilde{m}$.
    Let $s_i \in H^4(W_i; \mathbb{Q})$ be lifts of $ s \in H^4(M)$, such that $c_2(W_i) - m s_i \mapsto 0 \in H^4(M; \mathbb{Q})$.
    By exactness induced by 
    \begin{equation}
        M \rightarrow W_i \sqcup W_2 \rightarrow X_i
    \end{equation}
    there exists $ x_i \in H^4(X_i; \mathbb{Q})$, such that $x_i \mapsto (s_i, s_2) \in H^4(W_i; \mathbb{Q}) \oplus H^4(W_2; \mathbb{Q})$. 
    Then 
    \begin{equation}\begin{aligned}
        \int_{X_i} c_2(X_i)^2 = \int_{X_i} (c_2 (X_i) - mx_i)^2 = \int_{W_i} (c_2 (W_i) - ms_1)^2 + \int_{W_2} (c_2 (W_2) - m s_2)^2 ~\mathrm{mod} ~2\tilde{m}
    \end{aligned}\end{equation}
    The result then follows. 
\end{proof}

\subsection{Application to $\lieGroup{G}_2$-manifolds} 

In \cite{crowley2015new}, the authors define two invariants for 7 manifolds with $\lieGroup{G}_2$-structure, 
via a coboundary with $\lieGroup{Spin}(7)$-structure. 
Their methods are completely analogous to those presented here and are in fact the main inspiration. 

\begin{definition}\label{def_nu_prime} 
    Let $M$ be a closed 7 manifold with $\lieGroup{G}_2$-structure $\varphi$. 
    Let $W$ be a $\lieGroup{Spin}(7)$ coboundary to $M$. 
    Then 
    \begin{equation}
        \nu'(M, \varphi) := \chi(W) - 3 \sigma(W) ~~\mathSet{mod}~ 48 
    \end{equation}
\end{definition}

\begin{definition}\label{def_xi_prime} 
    Let $M$ be a closed 7 manifold with $\lieGroup{G}_2$-structure $\varphi$. 
    Let $p_M \in H^4(M)$ be the spin class, and set $m' \in \mathbb{N}$ be the greatest divisor of $p_M$ modulo torsion. 
    Define the set 
    \begin{equation}
        S_{m'} := \{ s \in H^4 (M) \colon p_M - m's \in TH^4(M) \} 
    \end{equation}

    Let $W$ be a $\lieGroup{Spin}(7)$ coboundary of $(M, \varphi)$ that is 3-connected. 
    We define a function $\xi' (\varphi): S_{m'} \rightarrow \mathbb{Q} / 3\tilde{m'} \mathbb{Z}$
    \begin{equation}
        \xi' (\varphi)(s) := 7 \chi(W) - \frac{45 \sigma(W)}{2} +  \int_W \frac{3 u_s ^2}{2} ~~\mathSet{mod} ~ 3 \tilde{m'} 
    \end{equation}
    where $ u_s  = p_W - m' \tilde{s} $, for some lift $\tilde{s} \in H^4(W) $ of $s$. 
\end{definition}

\begin{proposition}
\label{prop_compatible_g2_u3_structures}
    On a $7$-manifold $M$ with $\lieGroup{G}_2$-structure $\varphi$, there exists a compatible $\lieGroup{U}(3)$-structure $(v, g, \omega)$ with $c_1, c_3 =0$. 
    Moreover, 
    \begin{equation}\begin{aligned}
        \nu'(M, \varphi) &= \nu(M, v, g, \omega) & \xi(M, \varphi) &= \xi(M, v, g, \omega) 
    \end{aligned}\end{equation}
    In particular, this is independent of choice of compatible $\lieGroup{U}(3)$-structure. 
\end{proposition}

\begin{proof}
    A spin $7$-manifold admits an orthonormal pair of vector fields \cite[Theorem 1.1]{thomas1967postnikov}. 
    Thus a $\lieGroup{G}_2$-structure admits a reduction to an $\lieGroup{SU}(2)$-structure. 
    Fix such a reduction and extend trivially to a $\lieGroup{U}(3)$-structure. 
    Such a $\lieGroup{U}(3)$-structure clearly shares an $\lieGroup{SU}(3)$-reduction with $\lieGroup{G}_2$, 
    and $c_1, c_3 =0$. 
    
    In general $p_M = \frac{1}{2} (2c_2 + c_1 ^2)$, so for $M$, $p_M = c_2$ and thus $m' = m$.
    Consider an $\lieGroup{SU}(4)$ coboundary to the shared $\lieGroup{SU}(3)$ reduction. 
    The equality of the invariants is immediate by inspecting the remaining terms (note that $c_1(W) = 0 $ on an $\lieGroup{SU}(4)$-coboundary). 
\end{proof}

In light of this, we refer to a $\lieGroup{U}(4)$-coboundary of manifold with $\lieGroup{G}_2$-structure.
\begin{definition}\label{def_u4_cob_g2}
    Let $M$ be a closed $7$-manifold with a $\lieGroup{G}_2$-structure. 
    A $\lieGroup{U}(4)$-coboundary $W$ to $M$ is a coboundary such that the restriction of the $\lieGroup{U}(4)$-structure to $M$
    is a compatible $\lieGroup{U}(3)$-structure. 
\end{definition}

The following lemma is useful when checking compatible structures. 
\begin{lemma}
    \label{lem_compatible_g2_u3}
    Let $M$ be a 7-manifold with a $\lieGroup{G}_2$-structure $\varphi$ and $\lieGroup{U}(3)$-structure $(v, \omega, g)$ such that $c_1, c_3 =0$. 
    The structures are compatible provided that $\varphi \lrcorner v = \omega$, and $g_{\varphi} = g$. 
\end{lemma}

For a closed connected spin $7$-manifold $M$, let $\mathcal{G}_2$ denote the homotopy classes of $\lieGroup{G}_2$-structures on $M$. 
The homotopy classes of $\pi_0 \mathcal{G}_2 (M) \cong \mathbb{Z}$ (\cite[Lemma 1.1]{crowley2015new}). 
Denote the quotient by spin diffeomorphism by $\pi_0 \overline{\mathcal{G}}_2(M)$. 
There are various bounds available on the size $|\pi_0 \overline{\mathcal{G}}_2(M)|$ for various classes of $M$. 
We condense several into the following. 

\begin{proposition}\label{prop_nu_xi_complete_invariants}
    \cite[Theorem 1.12 \& 1.17]{crowley2015new}
    Let $M$ be a $2$-connected closed spin 7-manifold with torsion free cohomology such that $p_M \neq 0$, and $m = \mathrm{gd}(p_M)$.
    Then 
    \begin{equation}
        | \pi_0 \overline{\mathcal{G}}_2 (M) | = 24 \cdot \mathrm{Num}\left( \frac{m}{2^4 \cdot 7} \right) 
    \end{equation}
    and $(\nu, \xi)$ is a complete invariant of $\pi_0 \overline{\mathcal{G}}_2 (M)$ 
\end{proposition}

\section{The TCS construction} 
\label{sec_tcs_construction}

We restate some key definitions and results of the TCS of \cite{corti2015, kovalev2000twisted}. 

\subsection{TCS from ACyl CYs} 

Key to the construction is Calabi-Yau threefolds with tame asymptotic ends. 
We say $V_\infty ^{2n} $ is a Calabi-Yau (half)-cylinder if $V_\infty \cong \mathbb{R} ^+ \times X^{2n-1} $, 
is equipped with an $\mathbb{R}^+ $ invariant Calabi-Yau structure $(\omega_\infty, \Omega_\infty)$, 
such that induced metric $g_\infty$ is a product metric $dt^2 + g_X $ and $X$ is a smooth closed manifold.
$X$ is called the cross section of $V_\infty$.
The following is found in \cite[Definition 3.3]{corti2015}. 

\begin{definition}\label{def_accy} 
    Let $(V,\omega, \Omega)$ be a complete Calabi-Yau manifold. 
    Say that $V$ is an asymptotically cylindrical Calabi-Yau (ACyl CY) manifold if the following holds.
    There exists $(i)$ a compact set $K \subset V$, $(ii) $ a Calabi-Yau cylinder $V_\infty $, 
    and $(iii)$ a diffeomorphism $ \Phi: V_\infty \rightarrow V\setminus K $ such that for all $ k \in \mathbb{N}_0 $, and some $ \lambda >0 $ and as $t \rightarrow \infty$, 
    \begin{equation}\begin{aligned}
        \Phi^* \omega - \omega_\infty &= d \varrho,~~ \mbox{for $\varrho$ such that } | \nabla ^k \varrho | = O (e^{- \lambda t} ) \\
        \Phi^* \Omega - \Omega_\infty &= d \varsigma,~~ \mbox{for $\varsigma$ such that } | \nabla ^k \varsigma | = O (e^{- \lambda t} ) \\
    \end{aligned}\end{equation}
    where $\nabla$, and $ | \cdot | $ are defined in terms of $ g_\infty $ on $ V_\infty$. 
    We refer to $ V_\infty$ as the asymptotic end of $ V$. 
\end{definition}

The rate of convergence is important to the analysis justifying the existence of a torsion-free $\lieGroup{G}_2$-structure. 
We care only that the $\lieGroup{SU}(3)$-structure of the asymptotic end is an arbitrarily small perturbation from the cylindrical Calabi-Yau structure as we move along the neck.
In particular, the torsion free structure is homotopic to an $\lieGroup{SU}(3)$-structure that eventually agrees with that of the asymptotic end. 
For us the cross section is always $\Sigma \times S^1$, the product of a K3 and a circle. 
The asymptotic end has Calabi-Yau structure 
\begin{equation}\begin{aligned}\label{eq_accy_end}
    \omega_{\infty}  &= dt \wedge d \alpha + \omega ^I & \Omega_{\infty} &= (d\alpha - i dt) \wedge (\omega ^J + i \omega ^K)
\end{aligned}\end{equation}
for coordinates $(x, \alpha, t) \in \Sigma \times S^1 \times \mathbb{R}$, where $(\omega^I , \omega^J , \omega^K)$ is a hyper Kähler triple on $\Sigma$. 
This hyper Kähler K3 in the cross section of the asymptotic end is said to be the K3 at infinity of $V$. 

From an ACyl CY $V$ we can construct a torsion free $\lieGroup{G}_2$-structure on $M := V \times S^1$ given by 
\begin{equation}
    \varphi := d \beta \wedge \omega + \mathrm{Re}(\Omega) 
\end{equation}
where $\beta$ is the coordinate for the `external' circle factor. 
Sensibly extending our definitions of ACyl CY to $\lieGroup{G}_2$, the asymptotically cylindrical $\lieGroup{G}_2$ end of $M$ has $\lieGroup{G}_2$-structure
\begin{equation}\begin{aligned}\label{eq_acg2_end}
    \varphi_{\infty}&:= d \beta \wedge dt \wedge d \alpha + d \beta \wedge \omega_{} ^I +  d \alpha \wedge \omega_{} ^J +  d t \wedge \omega_{} ^K
\end{aligned}\end{equation}
Again, $\varphi$ is an arbitrarily small perturbation from $\varphi_\infty$ as we move along the neck.

We now turn our attention to defining a gluing on pairs of asymptotically cylindrical $\lieGroup{G}_2$ manifolds. 
\begin{definition}\label{def_hkr}
   Given hyper K\"ahler triples  $(\omega^I_\pm, \omega^J_\pm, \omega^K_\pm) $ on K3 surfaces $ \Sigma_\pm$ respectively, 
    a diffeomorphism $ r: \Sigma_+ \rightarrow \Sigma_- $ is called a Hyper K\"ahler Rotation (HKR) if $ r^* \omega^I _- = \omega^J _+$, $  r^* \omega^J _- = \omega^I _+$ and $ r^* \omega^K _- = -\omega^K _+$.
\end{definition}

Suppose that $M$ is an asymptotically cylindrical $\lieGroup{G}_2$-manifold obtained from ACyl CY $V$, and that $\Phi$ in Definition \ref{def_accy} has been specified. 
For fixed $T \gg 0$, we define forms $(\omega_T , \Omega_T )$ on $V$ 
\begin{equation}\begin{aligned}\label{eq_almost_su3_family} 
    \omega_{T}  & := \omega - d( \eta_T (t)  \varsigma) & \Omega_{T}  & := \Omega_\infty - d (\eta_T (t) \varrho)
\end{aligned}\end{equation}
where $\eta_T: \mathbb{R} \rightarrow [0,1]$ is a smooth cutoff function such that $\eta_T(t) = 0 $ for $t<T-1$, and $\eta_T(t) =1 $ for $t>T$. 
These forms are closed and interpolate between two torsion free $\lieGroup{SU}(3)$-structures on the neck. 
Note that $(\omega_{T} , \Omega_{T} )$ is not an $\lieGroup{SU}(3)$-structure. 
Despite this, for sufficiently large $T$, the $3$-form on $M$
\begin{equation}\label{eq_g2_structure_T} 
    \varphi_T := d \beta \wedge \omega_T  + \mathrm{Re}(\Omega_T) 
\end{equation}
is a $\lieGroup{G}_2$-structure on $M$ since the space of $3$-forms defining $\lieGroup{G}_2$-structures is open in the space of $3$-forms, 
and that $\varphi_T$ is a small perturbation from $\varphi$.  
The torsion of $\varphi_{T}$ is $O(e^{-\lambda T})$. 

Suppose we have a pair of ACyl CYs $V_\pm$ with cross sections $\Sigma_\pm \times S^1$, 
together with a HKR $r: \Sigma_+ \rightarrow \Sigma_-$ between K3s at infinity. 
We fix some coordinates $\Phi_\pm$ on the necks of $V_\pm$. 
Define a map 
\begin{equation}\begin{aligned}
    G_T: \Sigma_+ \times S^1 \times S^1 \times [T, T+1] & \rightarrow \Sigma_- \times S^1 \times S^1 \times [T, T+1]\\
        (x, \alpha, \beta, T+t) & \mapsto (r(x), \beta, \alpha, T+1-t) 
\end{aligned}\end{equation}

Note that $G_T ^* \varphi_\infty = \varphi_\infty$. 
Let $M(T)$ denote the truncation of $M$ at neck length $T+1$. 

\begin{definition}\label{def_tcs}
    For a pair of ACyl CY threefolds $V_\pm$ with cross section $\Sigma_\pm \times S^1$, together with a HKR $r: \Sigma_+ \rightarrow \Sigma_-$ between K3s at infinity, 
    their twisted connected sum $M := M_+ (T) \cup_{G_T} M_- (T) $ is a 7-manifold defined by the gluing of $M_\pm (T)$
    by the diffeomorphism $G_T$. 
\end{definition}
We endow $M$ with the closed $\lieGroup{G}_2$-structure $\varphi_T := \varphi_{+, T} \cup_{G_T} \varphi_{-, T}$.

\begin{theorem}\label{thm_tcs} 
\cite[Theorem 3.12]{corti2015} 
    Let $V_\pm$ be pair of ACyl CYs with a hyper Kähler rotation $r: \Sigma_+ \rightarrow \Sigma_-$ between K3s at infinity. 
    For sufficiently large $T$ there exists a torsion free perturbation of $\varphi_T$ within its cohomology class.
\end{theorem}

The resulting manifold $(M, \varphi)$ is closed with holonomy precisely $\lieGroup{G}_2$.

\subsection{ACyl CYs from algebraic geometry} 

ACyl CYs can be constructed from complex threefolds called building blocks. 

\begin{definition}\label{def_building_block}
   A building block is a non-singular algebraic threefold $Z$ together with a projective morphism
   $f: Z \rightarrow  \mathbb{P}^1 $ satisfying the following assumptions:
   \begin{enumerate}[topsep=0pt,itemsep=0pt,partopsep=0pt,parsep=0pt,label=(\roman*)]
      \item the anticanonical class $-K_Z  \in H^2 (Z) $ is primitive.
      \item $\Sigma := f^* (\infty ) $ is a smooth K3 with $ \Sigma \in |-K_Z | $.
   \end{enumerate}
   Identify $H^2 (\Sigma) $ with the K3 lattice $L$, ie choose a marking, 
   and let $N:= \mathrm{Im}(H^2(Z) \rightarrow H^2 (\Sigma) ) $.
   \begin{enumerate}[topsep=0pt,itemsep=0pt,partopsep=0pt,parsep=0pt,label=(\roman*)]\setcounter{enumi}{2}
      \item The inclusion $ N \hookrightarrow L $ is primitive.
      \item The group $H^3(Z) $ is torsion free.
   \end{enumerate}
    We may refer to a building block by $f:Z \rightarrow \mathbb{P}^1$, $(Z, \Sigma)$, or simply $Z$. 
\end{definition}

\begin{definition}\label{def_matching}
    Let $Z_\pm$ be a pair of complex threefolds, and $\Sigma_\pm \in |-K_{Z_\pm}|$ be smooth anticanonical divisors. 
    Let $k_\pm \in H^2(Z_\pm ; \mathbb{R})$ be Kähler classes and $\Pi_\pm \in H^2(\Sigma; \mathbb{R})$ be the 2-plane of type $(2,0) + (0,2)$.
    A diffeomorphism $ r: \Sigma_+ \rightarrow \Sigma_- $ is a \emph{matching}
    if $ r^* (k_-|_{\Sigma_-}) \in \Pi_+ $, $ k_+|_{\Sigma_+} \in r^*(\Pi_-) $, and $\Pi_+ \cap r^*(\Pi_-)$ is nonempty. 
\end{definition}
We refer to $r$ as a matching of $(Z_\pm, \Sigma_\pm, k_\pm)$ or a matching of $(Z_\pm, \Sigma_\pm)$ with respect to $k_\pm$. 
Where $k_\pm$ and/or $\Sigma_\pm$ exists but are not specified we refer to $r$ simply as a matching of $(Z_\pm, \Sigma_\pm)$ or of $Z_\pm$. 
The following result allows us to reformulate finding a HKR between building blocks to finding a matching. 

\begin{proposition}\label{prop_matching} 
\cite[Corollary 6.4]{corti2015}
    Let $(Z_\pm, \Sigma_\pm )$ be a pair of building blocks that have a matching  $r:\Sigma_+ \rightarrow \Sigma_-$ with respect to $k_\pm$. 
    Then $V_\pm := Z_\pm \setminus \Sigma_\pm $ admits a structure reduction to an ACyl CY structure (ie a torsion free $\lieGroup{SU}(3)$-structure)
    such that $r$ is a HKR on the K3s at infinity. 
\end{proposition}

We refer to the TCS of a pair of building blocks to mean the TCS of their associated ACyl CYs.

\section{A TCS coboundary} 
\label{sec_tcs_coboundary}

\subsection{A compatible $\lieGroup{U}(3)$-structure} 

We now give a description of a `nearby' $\lieGroup{G}_2$-structure on TCS $M$ which is homotopic to $\varphi$ given by Theorem \ref{thm_tcs}, 
and which aids our construction of describing a TCS $\lieGroup{U}(4)$-coboundary in the sense of Definition \ref{def_u4_cob_g2}. 

Suppose that $M$ is a TCS obtained from ACyl CYs $V_\pm$.
We endow $V_{\pm}$ with new $\lieGroup{SU}(3)$-structures $(\omega'_{\pm, T}, \Omega'_{\pm, T})$ such that 
\begin{equation}\begin{aligned}
    (\omega' _{\pm, T}, \Omega_{\pm, T} ') & = 
    \begin{cases} 
        (\omega_{\pm}, \Omega_{\pm}) & \mbox{ for } t< T-1 \\
        (\omega_{\pm, \infty}, \Omega_{\pm, \infty}) & \mbox{ for } t \geq T 
    \end{cases} 
\end{aligned}\end{equation}
Now we insist that $(\omega_{\pm, T} ', \Omega_{\pm,T} ')$ define $\lieGroup{SU}(3)$-structures, rather than remain closed 
in contrast to $(\omega_{\pm, T}, \Omega_{\pm, T})$ in Equation \ref{eq_almost_su3_family}. 
The space of $\lieGroup{SU}(3)$-structures above $ X = \Sigma \times S^1$ understood as reductions of frames of $TX \times \mathbb{R}$, forms a manifold.
This ensures the existence $(\omega' _{\pm, T}, \Omega_{\pm, T} ')$ for sufficiently large $T$. 

We endow $M_\pm = V_\pm \times S^1$ with the associated $\lieGroup{G}_2$-structures $\varphi_{\pm, T}'$. 
Define $\overline{V}_\pm$, to be the compact manifold with boundary obtained truncating $V_\pm$ at neck length $T$.
Define $\overline{M}_{\pm} := \overline{V}_{\pm} \times S^1$. 
We introduce $M_0 := \Sigma_+ \times T^2 \times [0, \tau]$ to be a piece of the neck and endow it with the asymptotically cylindrical $\lieGroup{G}_2$-structure $\varphi_\infty$. 
The parameter $\tau >0$ will be chosen later, and has no bearing on the homotopy class of the resulting $\lieGroup{G}_2$-structure. 
Note that $M_0$ has two boundary components, both diffeomorphic to $\Sigma_+ \times T^2$. 

We define gluings on the boundaries of $\overline{M}_\pm$ to $M_0$ by 
\begin{equation} \begin{aligned}
    G_\pm: \Sigma_\pm \times T^2 & \rightarrow \Sigma_+ \times T^2 \\ 
                     (x, \alpha, \beta) & \mapsto \begin{cases} 
                         (x, \alpha, \beta) & G_+ \\
                         (r^{-1}(x), \beta, \alpha) & G_- \\
                     \end{cases}
\end{aligned} \end{equation}
We note that in the next section we will construct a different hyper Kähler structure $M_0$, 
and with respect to which, $r$ will be an isomorphism of hyper Kähler manifolds. 
We think of $M$ as 
\begin{equation}\label{eq_M_glued_with_3_pieces} 
    M = \overline{M}_+ \cup_{G_+} M_0 \cup_{G_-} \overline{M}_-
\end{equation}
The $\lieGroup{G}_2$-structure $\varphi'$, formed from gluing together $\varphi_{\pm,T}'$ and $\varphi_\infty$, is homotopic to the torsion free $\lieGroup{G}_2$-structure $\varphi$.

\subsection{A TCS coboundary} 

Let $M$ be the TCS of ACyl CYs $(V_\pm, \Sigma_\pm)$ with HKR $r:\Sigma_+ \rightarrow \Sigma_-$. 
The idea behind the construction of a $\lieGroup{U}(4)$ coboundary is relatively straightforward---`rounding' the gluing of the coboundary used in \cite{crowley2014exotic} in order to handle a $\lieGroup{U}(4)$-structure.
We show that induced $\lieGroup{U}(3)$-structure on $M$ is compatible $\varphi'$. 

Let $W_\pm := \overline{V}_\pm \times D$, where $D \subset \mathbb{C}$ is the complex unit disc. 
We endow $W_\pm$ with the product $\lieGroup{SU}(4)$-structure 
determined by the $\lieGroup{SU}(3)$-structure $(\omega_{\pm, T}', \Omega_{\pm, T}')$  on $\overline{V}_\pm$, and the $SU(1)$-structure on $D$.

Note that $W_\pm$ are manifolds with corners - they are almost manifolds with boundary, but the boundary itself is the union of manifolds with boundary. 
We avoid needlessly engaging with the technicalities of the theory of manifolds with corners. 
We refer to the boundary components of $W_\pm$ as either \emph{internal} or \emph{external}. 
$\overline{V}_\pm \times S^1 \subset \partial W_\pm$ is the (external) boundary component identified with $\overline{M}_\pm$, 
while $\Sigma_\pm \times S^1 \times D \subset \partial W_\pm$ is the internal boundary. 
These two components meet along a common $6$-dimensional submanifold $\Sigma_\pm \times T^2$.
By construction the $\lieGroup{U}(4)$-structure on $W_\pm$ restricts to a $\lieGroup{U}(3)$-structure on $M_\pm$
compatible with $\varphi_\pm '$. 

We define a further manifold (with corners) $ W_0 := \Sigma \times Q$, where $\Sigma$ is a K3 and 
\begin{equation}\begin{aligned}
    Q &:=  \{ (z,w) \in \mathbb{C} ^2 \colon ~ |z|,|w|\leq 2, ~ (2 - |z|)^2 + (2 - |w|)^2 \geq 1 \} 
\end{aligned}\end{equation}
The boundary of $Q$ has three components 
\begin{equation}\begin{aligned}
    E_+ &:=  \{ (z,w) \in Q \colon ~ |z| = 2 \} \\
    E_- &:=  \{ (z,w) \in Q \colon ~ |w| = 2 \} \\
    Q_0 &:=  \{ (z,w) \in Q \colon ~ (2 - |z|)^2 + (2 - |w|)^2 = 1 \} 
\end{aligned}\end{equation}
We refer to $E_\pm$ as the internal boundary components while $Q_0$ is the (external) boundary of $W_0$. 
We will postpone a description of a $\lieGroup{U}(4)$-structure on $W_0$. 
It is clear that as smooth manifolds with boundary $M_0 \cong \Sigma \times Q_0$.  

Define the gluing maps between the internal boundaries of $W_\pm $ and $W_0$ 
\begin{equation} \begin{aligned}
    G_\pm: \Sigma_\pm \times S^1 \times D & \rightarrow \Sigma \times E_\pm \\ 
                     (x, \alpha, w) & \mapsto \begin{cases} 
                         (x, 2 e^{i \alpha}, w) & G_+ \\
                         (r^{-1}(x), w, 2 e^{i \alpha}) & G_- \\
                     \end{cases}
\end{aligned} \end{equation}
This is an extension of the gluing maps of $\overline{M}_\pm$ to $M_0$ to the interior of the internal boundaries. 

\begin{proposition}\label{prop_W_is_smooth_cob} 
    Let $W := W_+ \cup_{G_+} W_0 \cup _{G_-} W_-$ be the quotient space.
    Then $W$ is a smooth manifold with boundary and $\partial W \cong M$.
\end{proposition}

\begin{proof}
    The gluing map $G_+$ is a diffeomorphism onto its image.
    A neighbourhood of $\Sigma \times E_+ \subset W$ has a natural parameterization of $ \Sigma \times E_+ \times (-\varepsilon, \varepsilon)$,
    where $t \in (-\varepsilon, 0]$ belongs $W_+$, and otherwise $W_0$. 
    The analogous holds for $G_-$. 
    As the gluing is an extension of a gluing of $M$ in Equation \ref{eq_M_glued_with_3_pieces}, the boundary is smooth. 
    That the boundary is diffeomorphic to $M$ is clear from construction. 
\end{proof}

Let $M$ be the TCS of a pair of building blocks $Z_\pm$, 
and $W$ the coboundary given in Proposition \ref{prop_W_is_smooth_cob}. 
We can extend the embedding of $\overline{V}_\pm \rightarrow W$, given by $\overline{V}_\pm \rightarrow V_\pm \times \{0\} \subset W$, 
to the associated building blocks such that $Z_\pm \rightarrow W_\pm \cup_{G_\pm} W_0$.

\subsection{The $\lieGroup{U}(4)$-structure on the coboundary} 

One would like to define a $\lieGroup{U}(3)$-structure on $M$, use Lemma \ref{lem_compatible_g2_u3} to prove it is compatible with $\varphi'$, 
and prove it extends to a $\lieGroup{U}(4)$-structure on the interior of $W$. 
In practice, it seems easier to define a $\lieGroup{U}(4)$-structure on $W$ and check that the restriction to the boundary is compatible $\varphi'$.
The details are a little ugly. 

Suppose that $M$ is the TCS of building blocks $Z_\pm$, and $W$ is the coboundary of Proposition \ref{prop_W_is_smooth_cob}. 
We have an embedding of $\overline{V}_\pm \rightarrow W$ by identifying $\overline{V}_\pm \cong \overline{V}\pm \times \{0\} \subset W_\pm \subset W$. 
The complement of $\overline{V}_\pm$ in $Z_\pm$ is holomorphic to $\Sigma_\pm \times D$ where $D \subset \mathbb{C}$ is the complex unit disk. 
We define an embedding $j_\pm: Z_\pm \rightarrow W$ extending $\overline{V}_\pm$ by $\Sigma_\pm \times D \rightarrow \Sigma \times Q_0$.
The remainder of this Section provides the proof of the following proposition.

\begin{proposition}\label{prop_W_is_u4_cob}
    Let $W$ be as in Proposition \ref{prop_W_is_smooth_cob}.
    Then $W$ admits a $\lieGroup{U}(4)$-structure such that it is a $\lieGroup{U}(4)$-coboundary to $(M, \varphi')$. 

    Moreover, if the TCS is obtained from building blocks $Z_\pm$, then the pullback of the $\lieGroup{U}(4)$ structure via $j_\pm$
    reduces to a $\lieGroup{U}(3)$-structure that is homotopic to that induced by the complex structure on $Z_\pm$.  
\end{proposition}

We define a $\lieGroup{U}(4)$-structure on $W$ by considering each of its components in turn, 
and checking that they agree across the gluings. 
$W_\pm$ are each equipped with a product $\lieGroup{SU}(4)$-structure that on restricting to its external boundary agrees with the $\lieGroup{SU}(3)$-structure on $\overline{M}_\pm$. 

We define a $\lieGroup{U}(4)$-structure on $W_0$ while considering the constraints introduced 
by the gluing of the internal boundaries of $W_\pm$ and the $\lieGroup{G}_2$-structure on $M_0$. 
The $\lieGroup{U}(4)$-structure on $W_0$ will reduce to an $\lieGroup{SU}(2) \times \lieGroup{U}(2)$-structure and we consider $\Sigma$ and $Q$ in turn. 
$W_0$ can be viewed as a K3 fibration over base $Q$. 
The metric is the product metric $g_\Sigma + g_Q$ where $g_\Sigma$ is the Ricci flat metric shared by $\Sigma_\pm$;
and $g_Q$ is the metric on $Q$ that is yet to be determined.
Let $(\omega_I, \omega_J, \omega_K)$ be the hyper Kähler triple of $\Sigma_+$. 
We will endow each K3 fibre with a symplectic structure belonging to the $S^2$-family of forms this triple defines. 

Let $Q$ inherit the complex structure from $\mathbb{C}^2$, which will agree with the complex structure of the images of $G_\pm$ on $E_\pm$. 
Let $Q$ have hermitian metric $h_1 |dz|^2 + h_2 |dw|^2$, where $h_1$, and $h_2$ are positive real functions on $Q$. 
The constraints that this agrees with the images of $G_\pm$ are that $h_1  = 1/4 $, $h_2  = 1 $ on $E_+$ and $h_1  = 1 $, $h_2  = 1/4 $ on $E_-$. 

The precise description of the symplectic structure on K3 fibres, and on $h_i$ we finalize in the next point. 
That such functions exists is clear and any such choices will endow on $W_0$ a $\lieGroup{U}(4)$-structure.

As the $\lieGroup{U}(4)$-structures of $W_\pm$ agree with that of $W_0$ across their respective internal boundaries, 
so the structures can be glued to form $\lieGroup{U}(4)$-structures on $W$. 
As noted above, the $\lieGroup{U}(4)$-structure on $W_\pm$ restricts to $\lieGroup{U}(3)$-structures $(\partial_\beta, \omega_{\pm, T}', \Omega_{\pm,T}')$ on $\overline{M}_\pm$ respectively. 

It remains to specify $h_i$, and the symplectic structure of the K3 fibres of $W_0$, 
such that the restriction to $M_0$ is compatible with the $\lieGroup{G}_2$-structure $\varphi_{\infty}$. 
If $(v, \omega, g)$ is the $\lieGroup{U}(3)$-structure, then this amounts to checking that $\varphi_{\infty}\lrcorner v = \omega$, and $g_{\varphi_{\infty}} = g$ by Lemma \ref{lem_compatible_g2_u3}. 

Fix coordinates on $Q_0$,  
\begin{equation}\begin{aligned}
    f: S^1 \times S^1 \times [0, \textstyle \frac{\pi}{2}] & \rightarrow Q \\  
    ( \alpha, \beta, \vartheta) & \mapsto \left(( 2- \sin \vartheta) e^{i \alpha} , ( 2- \cos \vartheta) e^{i \beta}\right)
\end{aligned}\end{equation}
The derivative of $f$ is 
\begin{equation}
    df = \left( 
    \begin{array}{ccc}
        i(2 - \sin \vartheta) e^{i \alpha} & 0 & - \cos \vartheta e^{i \alpha} \\
        0 & i(2 - \cos \vartheta) e^{i \beta}  & \sin \vartheta e^{i \beta}  \\
    \end{array}
    \right) 
\end{equation}
Thus, the pullback metric on $Q_0$ by $f$ is
\begin{equation}
    g_{Q_0} := f^* g_{Q} = (h_1  \cos ^2 \vartheta + h_2  \sin ^2 \vartheta ) d \vartheta ^2 + h_1  (2 - \sin \vartheta)^ 2 d \alpha^2 + h_2   ( 2 - \cos \vartheta) ^2 d \beta ^2 
\end{equation}
The outward normal of $M_0$ is solely in the $Q$ component while the outward normal of $Q_0$ is  
\begin{equation}
    N = (h_2 \sin \vartheta e^{i \alpha}, h_1 \cos \vartheta e^{i \beta} ) ^T \in TQ 
\end{equation}
Thus $J(N)$ is a vector field on $Q_0$. 
Let $v'$ be the preimage under $df$ of $ J(N) $.
Then $v'$ is given in $(\vartheta, \alpha, \beta)$ coordinates by 
\begin{equation}
    v' = \frac{h_2 \sin \vartheta}{2-\sin \vartheta} \partial_\alpha  + \frac{h_1 \cos \vartheta}{2 - \cos \vartheta} \partial_\beta  
\end{equation}
We set the distinguished unit vector field $ v := v' / |v'| $ to be the normalized vector field.

Consider the symplectic form $ \omega_Q = g_Q (J_Q \cdot, \cdot)$ on $Q$ when pulled back to $Q_0$. 
\begin{equation}
     \omega_{Q_0} := f^* \omega_{Q} =  d \vartheta \wedge \left( - h_1  \cos \vartheta (2 - \sin \vartheta) d \alpha  + h_2  \sin \vartheta (2 - \cos \vartheta) d \beta \right)
\end{equation}
Thus the $\lieGroup{U}(3)$-structure on $M_0$ is $(v , g_{\Sigma} + g_{Q_0}, \omega_{\Sigma} + \omega_{Q_0}) $, where $\omega_\Sigma$ is the symplectic form on the K3 fibres. 

The $\lieGroup{G}_2$-form on $M_0$ is 
\begin{equation}
    \varphi_{\infty} = d \beta \wedge d t \wedge d \alpha + d \beta \wedge \omega_{} ^I +  d \alpha \wedge \omega_{} ^J +  d t \wedge \omega_{} ^K
\end{equation}
where $t=t(\vartheta)$ is a reparameterization of $\vartheta$. 
The metric associated to $\varphi_{\infty}$ is $ g_{\varphi_{\infty}} = g_{\Sigma_0} + dt^2 + d \alpha ^2 +  d \beta ^2 $. 
Thus, if this is to agree with $g_{\infty}$, then on $\infty$  
\begin{equation}\begin{aligned}
    h_1  & = ( 2 - \sin \vartheta)^{-2}, ~~h_2  = (2 - \cos \vartheta) ^{-2},  
\end{aligned}\end{equation}
while 
\begin{equation}\begin{aligned}
    \frac{d t}{d \vartheta} & = \left( \frac{\cos ^2 \vartheta}{(2- \sin \vartheta) ^2 } + \frac{\sin ^2 \vartheta}{(2-\cos \vartheta)^2} \right) ^{\frac{1}{2}} \\ 
\end{aligned}\end{equation}
As the righthand side is strictly positive, so $t(\vartheta)$ is strictly increasing. 

In terms of parameter $t$, the symplectic form and distinguished vector of the $\lieGroup{U}(3)$-structure on $M_0$ are 
\begin{equation}\begin{aligned}
    \omega & = \omega_{\Sigma} + \frac{d \vartheta}{dt} \left(\frac{\cos \vartheta }{2 - \sin \vartheta} d t  \wedge d \alpha +  \frac{\sin \vartheta }{2 - \cos \vartheta} d \beta \wedge d t  \right),  \\ 
    v & = \frac{d \vartheta}{d t} \left( \frac{\sin  \vartheta}{  2- \cos \vartheta} \partial_\alpha  + \frac{\cos \vartheta}{ 2 - \sin \vartheta} \partial_\beta \right) 
\end{aligned}\end{equation}
Contracting $\varphi_{\infty}$ with $v$ gives 
\begin{equation}
    \varphi_{\infty} \lrcorner v = \frac{d \vartheta}{d t} \left( \frac{\sin  \vartheta}{  2- \cos \vartheta} (   d \beta \wedge d t + \omega_2) + \frac{\cos \vartheta}{ 2 - \sin \vartheta} ( d t \wedge d \alpha + \omega_1) \right) 
\end{equation}

We now complete our description of the $\lieGroup{U}(4)$-structure on $W_0$.  
Let $\rho := \tan ^{-1} \left(\frac{2-|z|}{2-|w|} \right) $, and note that $\rho \circ f(\vartheta, \alpha, \beta) = \vartheta$. 
\begin{equation}
    \omega_{\Sigma} (\rho) := \left( \frac{\cos ^2 \rho}{(2- \sin \rho) ^2 } + \frac{\sin ^2 \rho}{(2-\cos \rho)^2} \right) ^{ -\frac{1}{2}}
\left( \frac{\sin  \rho}{  2- \cos \rho}  \omega_1 + \frac{\cos \rho}{ 2 - \sin \rho} \omega_2  \right) 
\end{equation}
and extend the definitions of $h_1$ and $h_2$ to the interior of $Q$ by 
\begin{equation}\begin{aligned}
    h_1  & = ( 2 - \sin \rho)^{-2}, ~~h_2  = (2 - \cos \rho) ^{-2} \\ 
\end{aligned}\end{equation}

Finally, the image of $Z_\pm$ restricted to $W_0$ is either $ \Sigma \times \{|z| \leq 2\} \times \{0\}$ of $\Sigma \times \{ 0 \} \times \{|w| \leq 2 \}$. 
It is clear then that the pullback of the $\lieGroup{U}(4)$-structure is homotopic to the $\lieGroup{U}(3)$-structure on $Z_\pm$ given by its complex structure. 
This concludes our proof of Proposition \ref{prop_W_is_u4_cob}.

\section{Topology of the TCS coboundary} 

To compute $\nu$ and $\xi$ for TCS manifolds we require an understanding of aspects of the cohomology. 
Let us fix the following notation.
\begin{definition}\label{def_lattice_notation}
    Let $L := H^2(\Sigma) $ be the abstract K3 lattice.
    Suppose we have a TCS from the building blocks $(Z_\pm, \Sigma_\pm)$
    Let $K_\pm := \mathSet{Ker}(H^2(Z_\pm ) \rightarrow L )$,
    and the polarization lattices $N_\pm := \mathSet{Im}(H^2(Z_\pm) \rightarrow L) $.

    For a pair of primitive embeddings $(N_\pm \rightarrow L)$ (which may or may not have come from a TCS)
    \begin{enumerate}
        \item $T_\pm := N_\pm ^\perp \subset L $, the transcendental lattices;
        \item $N_0 := N_+ \cap N_-$, the intersection lattice;
        \item $P := N_+ + N_- $ the span of the images of $N_\pm \hookrightarrow L$;
        \item $P_\pm := N_+ \cap N_- ^\perp$;
        \item $\Lambda_\pm := (N_+ ^\perp \cap N_-)^{\perp _P} = P_\mp ^{\perp _P} \subset P$.
    \end{enumerate}
\end{definition}
Note that our notation differs from \cite{corti2015} as here $K$ is the kernel of $H^2(Z) \rightarrow H^2(\Sigma)$, and not of $H^2(V)\rightarrow H^2(\Sigma)$. 

\subsection{Cohomology of $M$} 
\label{subsec_tcs_cohom}
A more detailed description of $H^*(M)$ is found in \cite[Theorem 4.8]{corti2015}. 
Applying Mayer-Vietoris to $(M_+, M_-)$ we get that $H^1(M) = 0 $, and  
\begin{equation}\begin{aligned}
    H^2 (M) & \cong  (K_+ / \mathSet{PD}(\Sigma_+)) \oplus (K_- / \mathSet{PD}(\Sigma_-)) \oplus N_0 \\
    H^3 (M) & \cong ( L/ P) \oplus (N_- \cap T_+ ) \oplus (N_+ \cap T_- ) \\
    & \qquad \oplus H^3(Z_+) \oplus H^3(Z_-) \oplus (K_+ \oplus K_- ) /\mathSet{PD}(\Sigma)    \\
    H^4 (M) & \cong (T_+ \cap T_-) \oplus (L/ (N_- + T_+))  \oplus (L/ (N_+ + T_-))\\ 
    & \qquad \oplus H^3(Z_+) \oplus H^3(Z_-) \oplus (K_+ ^ \vee \oplus_0 K_- ^\vee)    \\
\end{aligned}\end{equation}
where $  K_+ ^\vee \oplus_0 K_- ^\vee$ denotes the codimensional 1 subspace of $ (K_+ \oplus K_-)^\vee \cong K_+ ^\vee \oplus K_- ^\vee $, annhilating the subspace $\langle \mathrm{PD}(\Sigma) \rangle < K_+ \oplus K_-$.
In particular, the torsion part $TH^*(M)$ of $H^*(M)$ has the following form
\begin{equation}\begin{aligned}\label{eq_torsion_part_of_H_M} 
    TH^3(M) & \cong TH^3(Z_+) \oplus TH^3(Z_-) \oplus \mathrm{Tor}(L/P) \\
    TH^4(M) & \cong TH^3(Z_+) \oplus TH^3(Z_-) \oplus  \mathrm{Tor}(L/(N_+ + T_-) ) \oplus \mathrm{Tor}(L/(N_- + T_+)) 
\end{aligned}\end{equation}

\subsection{Cohomology of $W$} 
The long exact sequence of cohomology 
\begin{equation}\label{eq_rel_cohom} 
    H^*(W, M ;R) \rightarrow H^*(W;R) \rightarrow H^*(M;R) \rightarrow H^{*+1} ( W, M ;R) 
\end{equation}
relates the cohomologies of $W$ and $M$ and the relative cohomology of pair $(W, M)$. 

Let $ \tilde{W}_\pm := W_\pm \cup_{G_\pm } W_0$. 
We apply Mayer-Vietoris to $ W = \tilde{W}_+ \cup_{W_0} \tilde{W}_-$, noting that $ \tilde{W}_\pm \simeq Z_\pm $ and $W_0 \simeq \Sigma $.
Thus 
\begin{equation}
    H^*(W) \rightarrow H^*(Z_+) \oplus H^*(Z_-) \rightarrow H^*(\Sigma) \rightarrow H^{*+1} (W) 
\end{equation}
It follows that 
\begin{equation}\begin{aligned}
    H^1 (W) & \cong 0 & H^2 (W) & \cong K_+ \oplus K_- \oplus N_0  \\
    H^3 (W) & \cong L / (N_+ + N_-) \oplus H^3(Z_+) \oplus H^3(Z_-) & H^4 (W) & \cong H^4(Z_+) \oplus_0 H^4(Z_-)\\
\end{aligned}\end{equation}
where $  H^4(Z_+) \oplus_0 H^4(Z_-) $ denotes the codimensional 1 subspace of $ H^4(Z_+) \oplus H^4(Z_-) $
of pairs that share an image in $H^4(\Sigma) $. 

We recall that for CW complex pairs $(X,Y)$, $(A,C)$, $(B,D)$, suppose $X = A \cup B $, $Y = C \cup D$. 
Then a relative Mayer-Vietoris sequence is 
\begin{equation}
    H^*(X,Y) \rightarrow H^*(A,C) \oplus H^*(B,D) \rightarrow H^*(A \cap B, C \cap D)
\end{equation}
Let $\tilde{M}_\pm := M_\pm \cup M_0$. 
We apply the above sequence with $ A = W$, $B = W$, $C = W_+ \cup \tilde{M}_-$, and $D = W_- \cup \tilde{M}_+$. 
Then $A \cap B = W$ and $ C \cap D = M$. 
We have the following equivalences 
\begin{equation}\begin{aligned}
    (W, W_\mp \cup \tilde{M}_\pm) & \simeq (Z_\pm \times D,  Z_\pm \times S^1 ) \\ 
    (W , W_+ \cup \tilde{M}_- ) \cup (W , W_- \cup \tilde{M}_+) & \simeq (W , W_+ \cup W_- \cup M_0) \\ 
    (W , W_+ \cup W_- \cup M_0) & \simeq (\Sigma \times B^4, \Sigma \times S^3) \\ 
\end{aligned}\end{equation}
Thus, the Mayer Vietoris sequence becomes 
\begin{equation}
    H^{*-4} (\Sigma)  \rightarrow H^{*-2} (Z_+) \oplus H^{*-2} (Z_-) \rightarrow H^*(W, M) \rightarrow H^{*-3} (\Sigma) 
\end{equation}
This is in a sense dual to the decomposition used to compute $H^*(W)$.
We find that the long exact sequence gives the following exact sequences:
\begin{equation}\begin{aligned}
    & 0 \rightarrow H^0(Z_+) \oplus H^0 (Z_-) \rightarrow H^2 (W,M) \rightarrow 0 \\
    & 0 \rightarrow H^3(W,M) \rightarrow H^0(\Sigma) \rightarrow  H^2(Z_+) \oplus H^2(Z_-) \rightarrow H^4(W,M) \rightarrow 0 \\  
\end{aligned}\end{equation}
The map $H^0(\Sigma) \rightarrow  H^2(Z_+) \oplus H^2(Z_-)$ is generated by 
$\mathSet{PD}_{\Sigma} (\Sigma) \mapsto (\mathSet{PD}_{Z_+} (\Sigma) , \mathSet{PD}_{Z_-} (\Sigma) )$,
and in particular it is injective.
Hence
\begin{equation}\begin{aligned}
    H^1 (W,M) & \cong 0 & H^2 (W,M) &  \cong \left< \mathSet{PD}(Z_+),\mathSet{PD}(Z_-) \right>  \\
    H^3 (W,M) & \cong 0 &   H^4 (W,M) & \cong (H^2(Z_+) \oplus H^2(Z_-) ) / \left< \mathSet{PD}(\Sigma) \right> \\
\end{aligned}\end{equation}

From the dualities of the cohomology of $W$ and the relative cohomology of pair $(W,M)$, 
we have a complete description of their module structure.

\subsection{Characteristic classes} 

By Proposition \ref{prop_W_is_u4_cob} the Chern classes $c_j (\tilde{W}_\pm) \mapsto c_j (Z_\pm)$.
By $H^*(W) \rightarrow H^*(Z_+) \oplus H^*(Z_-)$ then $c_j(W) \mapsto c_j(Z_+) \oplus c_j(Z_-)$. 
The Euler characteristic $\chi(W) = \chi(Z_+) + \chi(Z_-) - 24$, since $\chi(\Sigma) = 24$.

\subsection{Products} 
To compute $\nu$ and $\xi$ we need to understand the following products: $c_1(W) c_3(W)$, $c_1(W)^2 c_2 (W)$, and $c_1 (W) ^4$. 
In addition, we need to make sense of the $c_2(W)^2$ term and compute the signature $\sigma(W)$.  
    
In the description of $H^*(W,M)$ we consider $(W, W_+ \cup \tilde{M}_-)$. 
We have a retract
\begin{equation}\begin{aligned}
    (Z_- \times D, Z_- \times S^1 ) & \rightarrow (W, W_+ \cup \tilde{M}_-)
\end{aligned}\end{equation}
while $(Z_- \times D, Z_- \times S^1)$ is topologically equivalent to the Thom space $Th(N_{W/Z_-})$ of the normal bundle of $Z_- \subset W$. 
We note that if $\iota_{\pm} : \Sigma \rightarrow Z_\pm$, then $\iota_{-} N_{W/Z_-} \cong N_{Z_+/ \Sigma}$. 

\begin{lemma}
    Let $A \subset B$ be a closed oriented codimension $k$ submanifold of closed oriented manifold $B$. 
    Then the map $H^*(A) \rightarrow H^{*-k}(B)$ given by  
    \begin{equation}
        H^*(A) \rightarrow H^{*-k}(Th(N_{B/A})) \rightarrow H^{*-k}(B,B-A) \rightarrow H^{*-k}(B)
    \end{equation}
    is adjoint to the restriction $H^*(B) \rightarrow H^*(A)$. 
\end{lemma}

It follows that 
\begin{center}
    \begin{tikzpicture}[scale=0.7]
        \tikzstyle{label} = []
        \node[label] (00) at (0,0)   {$ H^{*-2}(Z_\mp) $}; 
        \node[label] (01) at (6,0)   {$ H^*(W,M)$}; 
        \node[label] (02) at (12,0)  {$ H^*(W)$}; 
        \node[label] (10) at (0,-2)  {$ H^{*-2}(\Sigma)$}; 
        \node[label] (11) at (4,-2)  {$ H^{6-*}(\Sigma)^ \vee$}; 
        \node[label] (12) at (8,-2) {$ H^{6-*}(Z_\pm)^\vee$}; 
        \node[label] (13) at (12,-2) {$ H^*(Z_\pm)$}; 
        \draw[->, line width = 0.3mm] (00) -- (01) node[midway, above] {$ $};
        \draw[->, line width = 0.3mm] (01) -- (02) node[midway, above] {$ $};
        \draw[->, line width = 0.3mm] (10) -- (11) node[midway, above] {$ $};
        \draw[->, line width = 0.3mm] (11) -- (12) node[midway, above] {$ $};
        \draw[->, line width = 0.3mm] (12) -- (13) node[midway, above] {$ $};
        \draw[->, line width = 0.3mm] (00) -- (10) node[midway, rotate=90, above] {$ $};
        \draw[->, line width = 0.3mm] (02) -- (13) node[midway, rotate=90, above] {$ $};
    \end{tikzpicture}
\end{center}
By similar consideration we find that $H^{*-2}(Z_\pm) \rightarrow H^*(Z_\pm)$ factoring through $H^*(W,M) \rightarrow H^*(W)$ vanishes uniformly. 
This suffices to determine the products we require. 

We conclude that $H^2(W,M) \rightarrow H^2(W)$ is determined by $\mathrm{PD}(Z_\pm) \mapsto \mathrm{PD}(\Sigma_\mp)$. 
As $c_1(W) \cong c_1(Z_+) + c_1(Z_-)$ is in the image of $H^2(W,M) $, so it has a lift and so too $c_1(M) =0$. 
The products of $H^* _0(W)$ can be expressed in terms of the products of the building blocks. 
Thus
\begin{equation}
    c_1(W) ^2 \cong (\mathSet{PD}(\Sigma_+), \mathSet{PD}(\Sigma_-)) \in H^2(Z_+) \oplus H^2(Z_-) / \left< \mathSet{PD}(\Sigma) \right>  \cong H^4(W,M) 
\end{equation}
and so
\begin{equation}\begin{aligned}
    c_1(W) ^4 & = 0 & c_1(W) ^2 c_2(W) & = 48 & c_1(W) c_3(W) & = \chi(Z_-) + \chi(Z_+)  \\
\end{aligned}\end{equation}
Note that although $c_2(W) \notin H_0 ^4 (W)$, the triple product $c_1(W) (c_1(W) c_2(W)) \in H^8 (W,M)$ is still well defined.
For $H^6(W,M) \rightarrow H^6(W)$ we note only that it is onto, so that $c_3(W) \in H_0 ^6 (W)$. 

For $H^4(W,M) \rightarrow H^4(W)$, the map can be factored through the following diagram. 
\begin{center}
    \begin{tikzpicture}[scale=1]
        \tikzstyle{label} = []
        \node[label] (t1) at (0,0) {$(H^2(Z_+) \oplus H^2(Z_-))/ \left< \mathSet{PD}(\Sigma) \right> $}; 
        \node[label] (b1) at (0,-2) {$N_+ \oplus N_-$}; 
        \node[label] (b2) at (6,-2) {$ N_+ ^* \oplus N_- ^*$}; 
        \node[label] (t2) at (6,0) {$H^4(Z_+) \oplus_0 H^4(Z_-)$}; 
        \draw[->, line width = 0.3mm] (t1) -- (t2) node[midway, above] {$ $};
        \draw[->, line width = 0.3mm] (t1) -- (b1) node[midway, above] {$ $};
        \draw[->, line width = 0.3mm] (b1) -- (b2) node[midway, above] {$ \left( \begin{smallmatrix} 0 & \flat_- \\ \flat_+& 0 \end{smallmatrix} \right) $};
        \draw[->, line width = 0.3mm] (b2) -- (t2) node[midway, above] {$ $};
    \end{tikzpicture}
\end{center}
where $\flat_\pm = \iota_\mp ^* \circ Q_P \circ \iota_\pm$, for $Q_P$ the quadratic form on $P$. 
Thus the signature $\sigma(W) = 0$.
For computations it is also advantageous to note that 
letting $H_\pm := \mathrm{Im}(H^2(Z_\mp)\rightarrow H^2(Z_\pm))$, and $m_\pm := \mathSet{gd}(c_2(Z_\pm) ~\mathSet{mod}~(H_\pm))$,  
then $m := \mathrm{gd}(c_2(M)) = \mathrm{gcd}(m_+, m_-)$. 

In summary we have the following proposition. 
\begin{proposition}\label{prop_tcs_nu_xi}
    Let $(M, \varphi)$ be a TCS manifold and let $W$ be the TCS coboundary.
    Then $\sigma(W)=0$, $\chi(W) = \chi(Z_+) + \chi(Z_-) - 24$. 
    \begin{equation}\label{eq_nu_tcs_24}\begin{aligned}
        \nu(\varphi) = 24 ~(\mathSet{mod}~48)
    \end{aligned}\end{equation}
    If in addition $M$ has torsion free cohomology then 
    \begin{equation}\begin{aligned}
        \xi(\varphi) = \frac{3}{2} \int_W c_2 ^2 ~(\mathrm{mod}~3\tilde{m}) 
    \end{aligned}\end{equation}
\end{proposition}
Equation \ref{eq_nu_tcs_24} agrees with \cite[Theorem 1.7]{crowley2015new}.

\section{From semi-Fanos to TCS} 

We summarise some of the results in the literature of semi-Fanos we use, 
and their application to the properties of TCS manifolds. 

\subsection{Semi-Fano threefolds} 

\begin{definition}\label{def_weak_fanos}
    Let $Y$ be a smooth 3 dimensional complex algebraic variety with anticanonical class $-K_Y$.
    \begin{enumerate}
        \item If $-K_Y$ is ample, then $Y$ is \emph{Fano}. 
        \item If $-K_Y$ is big and nef (ie for all closed curves $C \subset Y$ $-K_Y \cdot C \geq 0$, and $(-K_Y)^3> 0$),
            then $Y$ is \emph{weak Fano}. 
        \item If $-K_Y$ is big and nef and semi-small (ie $|-K_Y|: Y \rightarrow \mathbb{P}^N$ at worst contracts finitely many divisors to curves)  
            then $Y$ is \emph{semi-Fano}. 
    \end{enumerate}
\end{definition}

All Fanos are semi-Fanos and all semi-Fanos are weak Fanos. 
Famously, there is a complete classification of Fano threefolds into 105 families,
but there are many deformation families of weak Fanos than Fanos.
Although we lack a classification of weak Fanos, many more examples are known. 
See authors \cite{jahnke2005threefolds, takeuchi2009weak, jahnke2011threefolds, arap2017existence} and references therein.

We collect some basic facts about weak Fanos. 
\begin{proposition}
    \label{prop_weak_fanos}
    Let $Y$ be a weak Fano threefold. 
    Then: 
    \begin{enumerate}
        \item \label{prop_weak_fanos_smooth_divisor} A general anticanonical divisor $\Sigma \in |-K_Y|$ is smooth. 
        \item \label{prop_weak_fanos_pic} $H^{0,i}(Y) = 0$ for $i>0$, and $\mathrm{Pic}(Y) \cong H^2(Y)$. 
        \item \label{prop_weak_fanos_c1_c2} $c_1(Y) \cdot c_2(Y) = 24$
        \item \label{prop_weak_fanos_degree} $h^0(Y, -K_Y) = g+2$ where $(-K_Y)^3 = 2g - 2$
        \item \label{prop_weak_ac_model} The anticanonical model $X$ of $Y$ is a Gorenstein Fano with at worst canonical singularities. 
        \item \label{prop_weak_transverse} $Y$ contains two smooth anticanonical divisors $\Sigma_i$ that intersect transversally 
            provided that its anticanonical model $X$ has very ample $-K_X$. 
        \item \label{prop_weak_fanos_canonical_curve} If $|-K_Y|$ contains two members $\Sigma_0, \Sigma_1$ that intersect transversally, 
            then curve $C = \Sigma_0 \cap \Sigma_1$ is a smooth canonically polarized curve (ie $K_C = -K_Y|_C$) 
            of genus $g$. 
        \item \label{prop_weak_fanos_semi_fano} $Y$ is a semi-Fano if and only if the bilinear form $(c_1(Y), \cdot, \cdot)$ on $H^2(Y)$ is non-degenerate. 
    \end{enumerate}
\end{proposition}

\begin{proof}
    For (\ref{prop_weak_fanos_smooth_divisor}) see \cite[Theorem 4.7]{corti2013asymptotically};
    for (\ref{prop_weak_fanos_pic}), (\ref{prop_weak_fanos_c1_c2}), (\ref{prop_weak_fanos_degree}) see \cite[Corollary 4.3]{corti2013asymptotically}; 
    for (\ref{prop_weak_ac_model}), (\ref{prop_weak_transverse}), (\ref{prop_weak_fanos_canonical_curve}) see \cite[Remark 4.10]{corti2013asymptotically}.
    For (\ref{prop_weak_fanos_semi_fano}), we note that if $Y$ is weak Fano and not semi-Fano there exists a divisor $D \subset Y$ such that $|-K_Y |$ contracts $D$ to a point. 
    Equivalently $c_1(Y) \cdot \mathrm{PD}(D) = 0 \in H^4(Y)$. 
\end{proof}

We call $g=Y(g)$ in the above Proposition \ref{prop_weak_fanos} the \emph{genus} of $Y$. 
It is equal to the genus of the canonically polarized curve $C$, and the genus of a smooth member $\Sigma \in |-K_Y|$. 
That is $g(C) = g(\Sigma) = g(Y)$. 

V. A. Iskovskikh \cite{iskovskih1977fano, iskovskih1978fano} lists the Fano threefolds for which the base locus of the anticanonical system is not empty. 
Jahnke-Radloff \cite[Theorem 1.1]{jahnke2006gorenstein} extend this to Gorenstein Fano threefolds with at worst canonical singularities.
Thus we can establish the base point freeness of the anticanonical model of a weak Fano threefold and from there discern the information of the weak Fano itself. 

\begin{proposition}
    Let $Y$ be the blowup of a rank 1 Fano $Y'$ of Fano index $r$ in a smooth embedded curve $C \subset Y'$ of degree $d$ and genus $g$. 
    Let $H$ be the pullback of the hyperplane class $H'$ on $Y'$ and $E$ the exceptional class. 
    $H$ and $E$ form a basis of the Picard lattice $N$, such that: 
    the product on $N$ is determined by $H^2 = \mathrm{deg}(Y')/ r^2$, $ E.H = d$ and $E^2 = 2g - 2$; 
    and the ample cone is spanned by rays $H$ and $G:= kH-E$ where $k$ is the least such integer such that $C$ is cut out by sections of $|kH|$. 
    In particular when $Y$ is known to be weak-Fano and not Fano then $G = -K_Y$.
\end{proposition}

See \cite[Lemma 4.5]{crowley2014exotic} and the succeeding remarks.  
We recall some results regarding the cohomology of blowups (see \cite[pg 605-611]{griffiths1978principles}).
Let $V \rightarrow X$ be a complex rank $r$ vector bundle over a smooth manifold $X$, and let $E = \mathbb{P}(V)$. 
The cohomology ring $H^*(E)$ is generated as an $H^*(X)$ algebra by 
\begin{equation}
    \left< \zeta \colon \textstyle \sum_{j=0} ^r ({-}1)^j c_j(V) \zeta^{r-j} = 0 \right> _{H^*(X)} 
\end{equation}
where $T \rightarrow \mathbb{P}(V)$ is the tautological line bundle, and $\zeta = c_1 (T)$. 
For a blowup $ \pi : Z \rightarrow Y$ along smooth complex submanifold $X \subset Y$ with exceptional divisor $E$, additively 
\begin{equation}\begin{aligned}
    H^*(Z) \cong \pi ^*H^*(Y) \oplus H^*(E) / \pi ^*H^*(C)
\end{aligned}\end{equation}
There is a natural isomorphism $E \cong \mathbb{P}(N_{X/Y})$ where $N_{X,Y}$ is the normal bundle of $X$ in $Y$. 
In the case of blowing up along a smooth complex curve $C$ in a threefold $Y$, $H^2(Z)$ is spanned by $H^2(Y)$ and $\zeta$. 
The product structure on even grades is determined by the product structure on $H^2(Y)$ and  
\begin{equation}\begin{aligned}
    \zeta^3 & = {-} \textstyle \int_C c_1(N_{C/Y}),&
    \pi^*(a) \cdot \zeta^2 & = {-} \textstyle \int_C a, & 
    \pi^*(b) \cdot \zeta & = 0  
\end{aligned}\end{equation}
for $a \in H^2(Y)$, $b \in H^4(Y)$. 
The Chern classes are given by 
\begin{equation}\begin{aligned}
    c_1(Z) & = \pi^*c_1(Y) - \zeta, &   
    c_2(Z) & = \pi^*(c_2(Y) + \mathrm{PD}(C) ) - \pi^* c_1(Y) \cdot \zeta 
\end{aligned}\end{equation}
Poincare duality implies that $H^2(Z) \cong H^4(Z) ^\vee$, and in computations it is helpful to express $H^4(Z)$ in the dual basis of $H^2(Z)$. 
In particular 
\begin{equation}\begin{aligned}
    \pi^*(a) \cdot c_2(Z) & = a \cdot c_2(Y) + \textstyle \int_C a, & 
    \zeta \cdot c_2(Z) & = \textstyle \int_C c_1(Y) 
\end{aligned}\end{equation}

In the case of Picard rank 1 Fano $Y$, $H^2(Y)$ is generated by the hyperplane class $H$. 
By definition, the (Fano) index $r$ of $Y$ is such that $c_1(Y) = r H$. 
As for any weak Fano threefold  $c_1(Y) \cdot c_2(Y) =24$, so $H \cdot c_2(Y) = 24/r$. 
For a smooth curve $C \subset Y$ of degree $d$ and genus $g$, the algebraic structure is determined by $H \cdot \zeta^2 = - d$ and $ \zeta^3 = - rd + \chi(C)$. 
The Chern classes are $c_1(Z) = rH - \zeta$, while $ H \cdot c_2(Z) = 24/r + d$ and $\zeta \cdot c_2(Z) = rd$.

\subsection{Building blocks of semi-Fano type} 

\begin{proposition}\label{prop_weak_bbs} 
    \cite[Proposition 5.7]{corti2013asymptotically}
    Let $Y$ be a weak Fano threefold, $C$ the base locus of a generic pencil in $|\Sigma_0 : \Sigma_1| \subset |-K_Y|$, 
    and assume that $C$ is smooth. 
    Let $Z$ be the blowup of $Y$ along $C$ and $f: Z \rightarrow \mathbb{P}^1$ the fibration induced by the pencil. 
    \begin{enumerate}
        \item The anticanonical class $-K_Z \in H^2(Z) $ is primitive. 
        \item The proper transform of the pencil $|\Sigma_0: \Sigma_1|$ is a fibration of K3s. 
        \item The restriction maps $H^2(Z) \rightarrow H^2(\Sigma)$ and $H^2(Y) \rightarrow H^2(\Sigma) $ have identical image. 
        \item $H^3(Z)$ is torsion free if and only if $H^3(Y)$ is torsion free. 
    \end{enumerate}
    Furthermore $\pi_1 (Z) = 0 $. 

    If, in addition, we suppose that $Y$ is a semi-Fano, then $H^2(Y) \rightarrow H^2(\Sigma)$ is a primitive embedding. 
    Thus, $Z$ is a building block. 
\end{proposition}

\begin{definition}
    A building block constructed from a Fano or semi-Fano as in Proposition \ref{prop_weak_bbs}
    will be said to be of Fano or semi-Fano type respectively. 
    A TCS of semi-Fanos will refer to the TCS of their associated building blocks and so on. 
\end{definition}

The even graded cohomology of a building block $Z$ obtained from semi-Fano $Y$ via Proposition \ref{prop_weak_bbs} is then determined as discussed above.
Let $C = \Sigma_0 \cap \Sigma_1$ be a canonical curve of genus $g$ to be blown up.
Let $\zeta \in H^2(Z)$ correspond to the exceptional class. 
Then $\zeta^3 = -2 K_Y ^3 = -2(2g-2)$ and $a \cdot \zeta^2 = - a \cdot c_1 ^2(Y)$. 
Note that $ a \cdot c_1 ^2 (Z) = a \cdot c_1^2 (Y) - \int_C a = 0$, 
and $\zeta \cdot c_1 ^2(Z) = \zeta^3 + 2 \int_C c_1(C) = 0$. 
Thus $c_1^2 (Z) =0$, which we could also see by noting that $Z$ is fibrated by anticanonical divisors. 
$c_2(Z)$ is as above.

\subsection{Matchings and configurations} 

A TCS manifold and coboundary are determined up to homotopy class of $G$-structures by a matching of building blocks. 
The question is then how does one find matchings---this is the matching problem. 
In \cite[Section 6]{corti2015}, the authors consider the `orthogonal' case (see discussion below). 
We require the more general `skew' case which has now been done in \cite{crowley2014exotic}. 

\begin{definition}
    \label{def_configuration} 
    Suppose $r: \Sigma_+ \rightarrow \Sigma_-$ is a matching (Definition \ref{def_matching}) between building blocks $(Z_\pm, \Sigma_\pm)$.
    After a choice of marking $H^2(\Sigma_+) \cong L$, we have a pair of primitive embeddings $N_\pm \rightarrow L$
    of the polarization lattices which we call the \emph{configuration} of $r$.
    It is well defined up to $ \lieGroup{O}(L)$. 
    
    A configuration is called \emph{orthogonal} if the reflections of $L(\mathbb{R})$ in $N_+$ and $N_-$ commute.
    If in addition $N_+ \cap N_-$ is trivial we call the configuration \emph{perpendicular}.
    If a configuration is not orthogonal, it is said to be \emph{skew}. 
\end{definition}

We consider the converse: 
given a pair of building blocks $(Z_\pm, \Sigma_\pm)$ with polarization lattices $N_\pm$,
which pairs of primitive embeddings $N_\pm \rightarrow L$ are realized as the configuration of some matching $r: \Sigma_+ \rightarrow \Sigma_-$. 
This is addressed in Proposition \ref{prop_genericity_match} for which we require some understanding of K3s.  

\begin{definition}
    The \emph{deformation family} $\mathcal{Y}(Y)$ of a weak Fano $Y$, with marking $i_Y : H^2(Y) \rightarrow N$ is a set of all $Y'$ 
    deformation equivalent to $Y$, and equipped with a compatible marking $i_{Y'} : H^2(Y') \rightarrow N$. 
    that is a surjective.
\end{definition}

As noted above the K3 fibres of a building block $Z$ of semi-Fano type are isomorphic to those in the pencil of the semi-Fano $Y$. 
In addition, the polarizations agree. 
By extension then we may refer to matchings of semi-Fanos to mean a matching of the associated building blocks.   
We endow $H^2(Y)$ with a lattice structure by taking the triple cup product and contracting with $-K_Y$. 
In this manner, for $\Sigma \in | -K_Y |$, $H^2(Y) \rightarrow \mathrm{Pic}(\Sigma)$ is a lattice embedding. 

We recall some terminology of lattice polarized K3s. 
The period domain of K3 surfaces is the space of oriented 2-planes in $L \otimes {\mathbb{R}}$. 
It can be identified with $ \{ \Pi \in \mathbb{P}(L \otimes \mathbb{C}) \colon \Pi ^2 = 0, \Pi \overline{\Pi} > 0 \} $, 
which inherits a complex structure. 
Given any $\Lambda \subset L$, the period domain of $\Lambda$-polarized K3s is 
$ D_\Lambda := \{ \Pi \in \mathbb{P}(\Lambda^\perp \otimes \mathbb{C} ) \colon \Pi ^2 = 0, \Pi \overline{\Pi} > 0 \} $. 

\begin{definition}\label{def_genericity}
    Let $N,\Lambda \subset L$ be primitive sublattices such that $N \subset \Lambda$. 
    Let $C \subset N(\mathbb{R})$ be a nonempty open subcone of the positive cone. 
    A set $\mathcal{Y}$ of $N$-marked threefolds is \emph{$(\Lambda, C)$-generic} if there exists a subset $U_\mathcal{Y} \subset D_\Lambda$
    that is the complement of a countable union of complex analytic submanifolds of positive codimension 
    with the property that: 
    For any $\Pi \in U$, $k \in C$, there is $Y \in \mathcal{Y}$, an AC divisor $\Sigma \subset Y$, 
    and a marking $h:H^2(\Sigma) \rightarrow L$ such that $\Pi$ is the period of $(\Sigma, h)$;
    the composition $H^2(Y) \rightarrow H^2(\Sigma) \rightarrow L$ equals a marking $i_Y: H^2(Y) \rightarrow N$;
    and $h^{-1}(k)$ is the image of the restriction to $\Sigma$ of a Kähler class on $Y$.
\end{definition}

Definition \ref{def_genericity} differs from \cite[Definition 6.17]{corti2015} where $U_\mathcal{Y}$ is a complement of a finite union of complex analytic submanifolds. 
The main advantage of this divergence is that the set $U_\Lambda \subset D_\Lambda$ for which there exists K3s with markings such that $\Lambda \cong \mathrm{Pic}(\Sigma)$
is generic in the sense used here, but it is not in the sense of \cite{corti2015}. 
By a \emph{genericity result} for a set $\mathcal{Y}$, 
we will mean a proposition of an exhaustible list of arithmetic conditions on a lattice $\Lambda$
that are sufficient to conclude that $\mathcal{Y}$ is $\Lambda$-generic. 

We note that Definition \ref{def_genericity} can be phrased (as in \cite[Definition 6.17]{corti2015}) in terms of families of building blocks, 
and so apply to cases where the building blocks are not of semi-Fano type. 
Likewise, the following result can be generalized to building blocks of other sources. 

\begin{proposition}\label{prop_genericity_match}
    Let $\mathcal{Y}_\pm$ be deformation families of semi-Fanos with polarizing lattices $N_\pm$, 
    and ample cones $C_\pm \subset N_\pm \otimes \mathbb{R}$. 
    Let $(N_\pm \hookrightarrow L)$ be primitive embeddings.
    With the notation of Definition \ref{def_lattice_notation}, assume that $ P_\pm \cap C_\pm $ is nonempty 
    and the $\mathcal{Y}_\pm$ are $(\Lambda_\pm, C_\pm)$-generic.

    Then there exists generic subsets $\mathcal{K}_\pm \subset P_\pm \cap C_\pm$ such that:
    for all pairs $k_\pm \in \mathcal{K}_\pm$, there exists $Y_\pm \in \mathcal{Y}_\pm$ 
    and smooth anticanonical divisors $\Sigma_\pm \rightarrow Y_\pm $, 
    such that there exists a matching of $(Z_\pm, \Sigma_\pm, k_\pm)$ with configuration $(N_\pm \rightarrow L)$.
\end{proposition}

The proof now appears in \cite[Proposition 5.8]{crowley2014exotic}.
The following is our first genericity result. 

\begin{proposition}
    \label{prop_beauville_genericity} 
    Let $Y$ be a semi-Fano threefold with Picard lattice $N$ and ample cone $C$.
    Then the deformation family $\mathcal{Y}(Y)$ is $(N, C)$-generic. 
\end{proposition}
See \cite[Theorem 6.8]{corti2013asymptotically} based on Beauville \cite{beauville2002fano}. 
In orthogonal configurations, the relevant $\Lambda$ is precisely isomorphic to the Picard lattice $N$. 
Thus we have sufficient genericity results for any orthogonal configuration. 
This has allowed for the mass production, particularly in the perpendicular cases (see \cite[Section 8]{corti2015}). 

Suppose that $P$ is isometric to the span of the images of $N_\pm \rightarrow L$ coming from some configuration. 
Then $P$ is a non-degenerate lattice of signature $(2, \mathrm{Rank}(P) -2)$, with primitive embeddings $N_\pm \rightarrow P$. 
It is easier, and sufficient for us to work with $P$ rather than embeddings into $L$, which is possible due to the following result by Nikulin. 
For a non-degenerate lattice $P$, we denote by $m(P)$ the minimal number of generators of $P^\vee /P$.  
\begin{theorem}\label{thm_nikulin_lattice_embedding}
\cite[Theorem 1.12.4 and Corollary 1.12.3]{nikulin1980integral}.
    Let $P$ be an even non-degenerate lattice of signature $(p_+, p_-)$, 
    and let $Q$ be an even unimodular lattice of indefinite signature $(q_+, q_-)$. 
    If $p_\pm \leq q_\pm $ and either 
    \begin{enumerate}
        \item $2 ~ \mathrm{Rank}(P) \leq \mathrm{Rank}(Q) $; or
        \item $\mathrm{Rank}(P) + m(P) < \mathrm{Rank}(Q) $;
    \end{enumerate}
    Then there exists a primitive embedding $P \hookrightarrow Q $. 
\end{theorem}

In the cases considered here $P$ is a non-degenerate lattice of rank $ \leq 4$ and so $P$ has signature $(2, \mathrm{Rank}(P)-2)$.
Thus it embeds primitively into the abstract K3 lattice $L$.
A primitive embedding of $P \hookrightarrow L$ induces primitive $N_\pm \hookrightarrow L$. 
Note that the converse is not necessarily true. 
That is, we may find overlattice refinements $P \hookrightarrow \tilde{P}$, where $|\tilde{P}/P|$ is finite and non-zero, 
and then primitively embed $\tilde{P} \hookrightarrow L$. 
This leads to TCS manifolds with torsion in $H^3$, and we shall not consider this further. 

In the case of perpendicular, or orthogonal matchings, $P$ is the pushout of $N_\pm$. 
It follows that perpendicular pushouts involving Picard lattices of rank $\leq 2$ always correspond to a matching.
By specifying shared isometric primitive sublattice, $N_0 \rightarrow N_\pm$ we can define an orthogonal pushout $P = (N_+ \oplus N_-)/N_0$. 
One must have selected $N_0$ to be orthogonal to some region of the relevant cone in each $N_\pm$, 
and check that $P$ is integrable. 
If met, then it corresponds to a matching. 

A classification of closed smooth 7-manifolds with non-zero $b_2$ is under investigation, and off the shelf results are currently insufficient.
As $b_2(M) \geq  \mathrm{Rank}(N_0)$, our examples all have trivial intersection lattice.
It is quite conceivable that an enhanced classification of closed smooth 7-manifolds would lead to orthogonal matching examples demonstrating the use of $\xi$, 
which would avoid more involved genericity results. 

For the skew case with trivial intersection we start with taking $P := \left(\begin{smallmatrix} N_+ & D  \\ D^T & N_- \end{smallmatrix}\right)$ 
as a block matrix and aim to find the solutions of matrix $D$ such that $P$ is a nondegenerate lattice of rank $p = \mathrm{Rank}(N_+) + \mathrm{Rank}(N_-)$ and signature $(2, p - 2)$. 
In addition, we also check that $W_\pm \cap C_\pm$ is nonempty. 
For each solution, we have a $\Lambda_\pm$ for which we require a genericity result for the deformation family $\mathcal{Y}$. 
The integrality constraint ($D_{ij} \in \mathbb{Z}$) together with $W_\pm \cap C_\pm \neq \varnothing$
determine that the number of solutions is always finite.

\section{Examples} 
\label{sec_examples}

\subsection{Some genericity results} 

For our examples we have skew matchings that require genericity results which do not follow from Beauville. 
For a Picard rank 1 Fano threefold $Y'$, let $H$ denote the fundamental class, ie the hyperplane class. 
Let $Y$ be a Picard rank 2 Fano or semi-Fano threefold. 
We denote its Picard lattice by $N$. 
In the case that $Y$ is obtained by blowing up a Picard rank 1 Fano $Y'$,
we denote by $H \in N$ the pullback of the fundamental class of $Y'$; $E \in N$ the exceptional class; 
and primitive class $G \in N$ such that $G,H$ span the ample cone. 
In all cases considered, $G,H$ are a basis of $N$. 

The matchings involve the following semi-Fano families. 
\begin{enumerate}[topsep=0pt,itemsep=0pt,partopsep=0pt,parsep=0pt,label=(Y\arabic*)]
    \item \label{example_prebb_1_Q_0_6} The blowup of a smooth quadric in a smooth rational curve of degree 6 \cite[No. 109]{cutrone2013towards}.
        Here $G := 3 H - E$ and $N  = \left(\begin{smallmatrix} 16 & 12  \\ 12 & 6 \end{smallmatrix}\right)$. 
        \item \label{example_prebb_1_PP3_8_2} The blowup of $\mathbb{P}^3$ in a smooth curve of degree 8 and genus 2 \cite[No. 49]{cutrone2013towards}.
        Here $G := 4 H - E$ and $N  = \left(\begin{smallmatrix} 2 & 8  \\ 8 & 4 \end{smallmatrix}\right)$. 
        \item \label{example_prebb_1_PP3_11_14} The blowup of $\mathbb{P}^3$ in a smooth curve of degree 11 and genus 14 \cite[No. 52]{cutrone2013towards}.
        Here $G := 4 H - E$ and $N  = \left(\begin{smallmatrix} 2 & 5  \\ 5 & 4 \end{smallmatrix}\right)$. 
        \item \label{example_prebb_1_PP1_x_PP2} The Fano threefold $\mathbb{P}^1 \times \mathbb{P}^2$ \cite[Ch.12, Table 2, row 34]{parshin1999algebraic}.
        Let $G$ be the pullback of the hyperplane class on $\mathbb{P}^1$ and $H$ be the pullback of the hyperplane class on $\mathbb{P}^2$. 
        With respect to basis $(G,H)$, $N  = \left(\begin{smallmatrix} 0 & 3  \\ 3 & 2 \end{smallmatrix}\right)$. 
        and $G,H$ form a basis of the ample cone. 
        \item \label{example_prebb_1_10} The Fano threefold $Y$ of Picard rank 1, of first species, and genus 10 \cite[Ch.12, Table 1, row 9]{parshin1999algebraic} 
            It has degree 18.
\end{enumerate}
All but \ref{example_prebb_1_10} require some genericity result to justify the configurations used in our examples.

Our proofs of the genericity results use critically the description of the family of the relevant semi-Fano. 
Thus we do not get such a general one-size-fits-all result of Beauville. 
However, similarities in descriptions of deformation families lead to similarieties in genericity results.
For example, in the case that $Y$ is the blowup of some $Y'$ in a smooth curve of set degree and genus, 
a genericity result can be obtained by first embedding a lattice polarized K3 into $Y'$, 
then checking there exists a smooth curve representative of some class that has the correct degree and genus. 
\cite[Section 12]{parshin1999algebraic} contains a tabulated summary of the classification of Fano threefolds. 
Of the 36 deformation families of rank 2 Fanos, 27 are presented as blowups along curves embedded in rank 1 Fanos. 
Moreover, the presentation of deformation families of rank 2 weak Fanos are also in this form. 
Of the remaining nine families of rank 2 Fano threefolds they are either: projective bundles over $\mathbb{P}^2$; 
divisors of $\mathbb{P}^2 \times \mathbb{P}^2$; or double covers of other rank 2 Fanos. 
Thus one suspects that a systematic approach can be conducted if desired. 

We first state some useful results concerning K3s. 
The following is based on \cite[Chapter 3]{reid1996chapters}. 

\begin{proposition}\label{prop_reid_k3s}
    Let $\Sigma$ be a $\Lambda$-polarized K3, and suppose $H \in \Lambda$ is nef and big. 
    Then either:
    \begin{enumerate}
        \item $|H|$ is base point free, or
        \item $H$ is monogonal, ie $H = aD + \Gamma$ where $a \geq 1$, $D^2 = 0$, $D.\Gamma = 1$ and $\Gamma^2 = -2$. 
    \end{enumerate}

    In the case that $H$ is not monogonal, then either 
    \begin{enumerate}
        \item the morphism enduced by $|H|$ birational onto its image and an isomorphism away from a finite union of $(-2)$-curves, or 
    \item $|H|$ is hyper-elliptic, ie one of the following holds 
        \begin{enumerate}
            \item $H^2 = 2 $ and $\Sigma$ is a double cover of $\mathbb{P}^2$, or
            \item $H = 2B $ with $B^2 = 2$ and $\Sigma$ is a double cover of the veronese surface, or
            \item $\Sigma$ has an elliptic pencil $|E|$ with $H.E = 2$.
        \end{enumerate}
    \end{enumerate}
\end{proposition}

\begin{corollary}\label{cor_veryample_k3}
    Let $\Sigma$ be a $\Lambda$-polarized K3, and suppose $H \in \Lambda$ is nef and big. 
    If 
    \begin{enumerate}
        \item \label{cor_vamp_Hgeq4} $H^2 \geq 4$ 
        \item \label{cor_vamp_prim} $H$ is primitive
        \item \label{cor_vamp_sq0} $\nexists D \in \Lambda$ such that $H.D = 2 $ and $D^2 = 0$
        \item \label{cor_vamp_sqn2} $\nexists D \in \Lambda$ such that $H.D = 0 $ and $D^2 = -2$
    \end{enumerate}
    Then $H$ is very ample.
\end{corollary}

\begin{proof}
    If $H$ were monogonal then there exists a class $E$ such that 
    $H.E = 1$ and $E^2 = 0$, but then class $2E$ contradicts (\ref{cor_vamp_sq0}). 
    Thus $H$ is base point free.
    By (\ref{cor_vamp_Hgeq4}), $|H|$ does not enduce a double cover of $\mathbb{P}^2 $; 
    by (\ref{cor_vamp_prim}), $|H|$ does not enduce a double cover of the Veronese surface; 
    and by (\ref{cor_vamp_sq0}), $\Sigma$ cannot contain an elliptic pencil $|E|$ with $H.E = 2$. 
    Thus $H$ is not hyper-elliptic. 
    Finally, (\ref{cor_vamp_sqn2}) rules out the possibility that $|H|$ contracts any $(-2)$-curves. 
    Thus $H$ is very ample.
\end{proof}

\begin{proposition}\label{prop_smooth_representative}
    Let $\Sigma$ be a $\Lambda$-polarised K3 with very ample class $H$.
    Suppose $E \in \Lambda$ is such that $E.H>0$. 
    \begin{enumerate}
        \item If $E^2 = -2$ then $E$ irreducible; or 
        \item If $E^2 = 0 $ then $\nexists D \in \Lambda $ such that $D^2 = -2 $, $0< D.H< H.E$, and $ E.D < 0 $; or 
        \item If $E^2 > 0 $ then $\nexists D \in \Lambda $ such that $D^2 = -2 $, and either $E.D < 0$; or there exists an $a \geq 2$ 
            and class $E'$ such that $E = aE' + D $, $(E')^2 = 0 $ and $ E'.D = 1$, 
    \end{enumerate}
    then $E$ is represented by a smooth curve $C$ of degree $d$.
\end{proposition}

\begin{proof}
    By construction $E.H > 0 $, so in each case \cite[Section 3.7]{reid1996chapters} implies $h^0 (D) > 0 $. 
    In particular $E$ is effective. 
    Also we have that any irreducible class $D \in \Lambda$ with $D^2 = -2 $ is represented by a smooth rational curve. 
    Now consider each case 
    \begin{enumerate}
        \item By \cite[section 3.8 (a)]{reid1996chapters}  we can decompose any effective class as $E = M + F$, 
            where $M$ is effective and nef with $M^2 \geq D^2 $ and 
            $F = \sum_i n_i \Gamma_i $ an effective sum of $(-2)$-curves. 
            Thus $H.E = H.M + \sum_i n_i H. \Gamma_i $. 
            In excluding this possibility, $E$ is irreducible. 
        \item By \cite[section 3.8 (b)]{reid1996chapters} provided that $E$ is nef then $E$ has a base point free complete linear system. 
            It would fail to be nef if and only if there exists a $(-2)$-curve on which $E$ were negative. 
            The stipulation rules this out.
        \item By \cite[section 3.8 (d)]{reid1996chapters} provided that $E$ is nef and not monogonal then $E$ has a base point free complete linear system.  
            Again we exclude precisely these possibilities by our stipulation. 
    \end{enumerate}
    In the second two cases, we apply Bertini's theorem that says a base point free complete linear system has a smooth divisor. 
\end{proof}

In addition, we have the following result. 
\begin{lemma}\label{cor_removing_n2curves}
    Let $\Sigma$ be a K3 surface, $D$ an effective class, and $\Gamma$ a $(-2)$-curve. 
    If $D . \Gamma < 0 $ then $D - \Gamma$ is effective. 
\end{lemma}
We now present genericity results relevant to our needs. 

\begin{proposition}\label{prop_genericity_wf_1_16_6_0}
    Let $Y$ be a semi-Fano obtained by blowing up a quadric in $\mathbb{P}^4$ along a smooth rational curve of degree $6$. 
    Let $\mathcal{Y} = \mathcal{Y}(Y)$ be the deformation family associated to $Y$. 
    Let $\Lambda$ be an overlattice of $N \rightarrow \Lambda$ in which $N$ is primitive. 
    If all:
    \begin{enumerate}
        \item \label{prop_genericity_wf_1_16_6_0_i} $\nexists D \in \Lambda$ such that $H.D = 2 $ and $D^2 = 0$
        \item \label{prop_genericity_wf_1_16_6_0_ii} $\nexists D \in \Lambda$ such that $H.D = 0 $ and $D^2 = -2$
        \item \label{prop_genericity_wf_1_16_6_0_iii} $\nexists D \in \Lambda$ such that $H.D = 3 $ and $D^2 = 0$
        \item \label{prop_genericity_wf_1_16_6_0_iv} $\nexists D \in \Lambda$ such that $0< H.D < 6$ and $D^2 = -2$, and $E.D < 0$. 
    \end{enumerate}
    hold, then $\mathcal{Y}$ is $\Lambda$-generic. 
\end{proposition}

\begin{proof}
    By assumption that $N \hookrightarrow \Lambda$ is a primitive embedding, so $H \in \Lambda$ is again primitive. 
    As $H^2 = 6 \geq 4$, together with hypotheses (\ref{prop_genericity_wf_1_16_6_0_i}) and (\ref{prop_genericity_wf_1_16_6_0_ii})
    we conclude that $H$ is very ample by Corollary \ref{cor_veryample_k3}. 
    Thus $|H|: \Sigma \rightarrow \mathbb{P}^4$ is a smooth embedding. 
    By \cite[Lemma 2.4]{fukuoka2017existence}, hypothesis (\ref{prop_genericity_wf_1_16_6_0_iii}) ensures that $|H|(\Sigma)$ is an anticanonical divisor of a smooth quadric $Q \subset \mathbb{P}^4$. 
    Hypothesis (\ref{prop_genericity_wf_1_16_6_0_iv}) ensures that $E$ is an irreducible class and so represented by a smooth $(-2)$-curve $C$ by Proposition \ref{prop_smooth_representative}.
    Blowing up $Q$ along $C$ gives us a member of $\mathcal{Y}$, while the proper transform of $|H|(\Sigma)$ is an anticanonical divisor that is again isomorphic to $\Sigma$. 
\end{proof}

\begin{proposition}\label{prop_genericity_wf_1_17_8_2}
    Let $Y$ be a semi-Fano obtained by blowing up $\mathbb{P}^3$ along a smooth curve of degree $8$ and genus $2$. 
    Let $\mathcal{Y} = \mathcal{Y}(Y)$ be the deformation family associated to $Y$. 
    Let $\Lambda$ be an overlattice of $N \rightarrow \Lambda$ in which $N$ is primitive. 
    If all:
    \begin{enumerate}
        \item \label{prop_genericity_wf_1_17_8_2_i} $\nexists D \in \Lambda$ such that $H.D = 2 $ and $D^2 = 0$
        \item \label{prop_genericity_wf_1_17_8_2_ii} $\nexists D \in \Lambda$ such that $H.D = 0 $ and $D^2 = -2$
        \item \label{prop_genericity_wf_1_17_8_2_iii} $\nexists D \in \Lambda$ such that $E.D <  $ and $D^2 = -2$
        \item \label{prop_genericity_wf_1_17_8_2_iv} $\nexists D \in \Lambda$, such that $D^2 = -2$, $(E - D)^2 = 0$, and $E.D =\frac{1}{2} E^2 -1$.
    \end{enumerate}
    hold, then $\mathcal{Y}$ is $\Lambda$-generic. 
\end{proposition}

\begin{proof}
    By assumption that $N \hookrightarrow \Lambda$ is a primitive embedding, so $H \in \Lambda$ is again primitive. 
    As $H^2 = 4$, together with hypotheses (\ref{prop_genericity_wf_1_17_8_2_i}) and (\ref{prop_genericity_wf_1_17_8_2_ii})
    we conclude that $H$ is very ample by Corollary \ref{cor_veryample_k3}. 
    Thus $|H|: \Sigma \rightarrow \mathbb{P}^3$ is a smooth embedding. 
    Hypothesis (\ref{prop_genericity_wf_1_17_8_2_iii}) ensures that $E$ is nef. 
    Hypothesis (\ref{prop_genericity_wf_1_17_8_2_iv}) ensures that $E$ is not monogonal.
    Thus by Proposition \ref{prop_smooth_representative}, $E$ is represented by a smooth $C$. 
    This curve has degree $E.H = 8$ and genus $\frac{1}{2}E^2 + 1 = 2$. 
    Blowing up $\mathbb{P}^3$ along $C$ gives us a member of $\mathcal{Y}$, while the proper transform of $|H|(\Sigma)$ is an anticanonical divisor that is again isomorphic to $\Sigma$. 
\end{proof}

\begin{proposition}\label{prop_genericity_wf_1_17_11_14}
    Let $Y$ be a semi-Fano obtained by blowing up $\mathbb{P}^3$ along a smooth curve of degree $8$ and genus $2$. 
    Let $\mathcal{Y} = \mathcal{Y}(Y)$ be the deformation family associated to $Y$. 
    Let $\Lambda$ be an overlattice of $N \rightarrow \Lambda$ in which $N$ is primitive. 
    If all:
    \begin{enumerate}
        \item \label{prop_genericity_wf_1_17_11_14_i} $\nexists D \in \Lambda$ such that $H.D = 2 $ and $D^2 = 0$
        \item \label{prop_genericity_wf_1_17_11_14_ii} $\nexists D \in \Lambda$ such that $H.D = 0 $ and $D^2 = -2$
        \item \label{prop_genericity_wf_1_17_11_14_iii} $\nexists D \in \Lambda$ such that $E.D <  $ and $D^2 = -2$
        \item \label{prop_genericity_wf_1_17_11_14_iv} $\nexists D \in \Lambda$, such that $D^2 = -2$, $(E - D)^2 = 0$, and $E.D =\frac{1}{2} E^2 -1$.
    \end{enumerate}
    hold, then $\mathcal{Y}$ is $\Lambda$-generic. 
\end{proposition}

\begin{proof}
    Completely analogous to Proposition \ref{prop_genericity_wf_1_17_8_2}.  
\end{proof}

Finally, we require a genericity result for $\mathbb{P}^1 \times \mathbb{P}^2$. 
\begin{proposition}\label{prop_genericity_m_2_34}
    Let $Y = \mathbb{P}^1 \times \mathbb{P}^2$. 
    The associated deformation family of $Y$ is simply $\{Y \}$.
    Let $\Lambda$ be an overlattice of $N \rightarrow \Lambda$ in which $N$ is primitive. 
    If all:
    \begin{enumerate}
        \item \label{prop_genericity_m_2_34_i} $\nexists D \in \Lambda$ such that $(G+H).D = 2 $ and $D^2 = 0$
        \item \label{prop_genericity_m_2_34_ii} $\nexists D \in \Lambda$ such that $(G+H).D = 0 $ and $D^2 = -2$
        \item \label{prop_genericity_m_2_34_iii} $\nexists D \in \Lambda$ such that $0< (G+H).D < 3 $ and $D^2 = -2$, and $G.D < 0$. 
        \item \label{prop_genericity_m_2_34_iv} $\nexists D \in \Lambda$ such that $0< (G+H).D < 5 $ and $D^2 = -2$, and $H.D < 0$. 
    \end{enumerate}
    hold, then $\mathcal{Y}$ is $\Lambda$-generic.  
\end{proposition}

\begin{proof}
    As $(G+H)^2 = 8$, $G+H$ is primitive then by the hypotheses (\ref{prop_genericity_m_2_34_i}) and (\ref{prop_genericity_m_2_34_ii}), 
    $G+H$ is very ample by Corollary \ref{cor_veryample_k3}.
    Then $G,H$ are both effective since their squares $\geq -2$ and their products with $(G+H)$ are each positive. 

    By (\ref{prop_genericity_m_2_34_iii}) $G$ is nef.
    If it were not nef, there would exist an effective $(-2)$-curve $\Gamma$ such that $G. \Gamma <0 $. 
    By Corollary \ref{cor_removing_n2curves}, $G-\Gamma$ is effective. 
    So 
    \begin{equation}
        0< (G+H).(G - \Gamma) = 3 - (G+H).\Gamma
    \end{equation}
    which contradicts (\ref{prop_genericity_m_2_34_iii}). 
    Similarly (\ref{prop_genericity_m_2_34_iv}) implies that $H$ is nef. 

    As $G$ is nef and $G^2 = 0$, by \cite[Section 3.8]{reid1996chapters} there exists an effective class $E$ 
    such that $G = aE$ and $|E|$ is a free pencil. 
    As $G$ is primitive, $a=1$, and so $h^0(G) = 2$. 
    As $H$ is nef and big, \cite[Section 3.8]{reid1996chapters} says that $h^0(H) = \frac{1}{2} H^2 +2 =3$. 
    Likewise $h^0(G+H) = 6 $. 

    We have the following diagram. 
    \begin{center} 
        \begin{tikzpicture}[scale=1]
            \tikzstyle{label} = []
            \node[label] (1) at (0,0) {$\Sigma$}; 
            \node[label] (21) at (3,1) {$\mathbb{P}^1 \times \mathbb{P}^2 $}; 
            \node[label] (22) at (3,-1) {$\mathbb{P}^5$}; 
            \draw[->, line width = 0.3mm] (1) -- (21) node[midway, above, rotate=18] {$(|G|, |H|)$};
            \draw[->, line width = 0.3mm] (1) -- (22) node[midway, below, rotate=-18] {$|G + H|$};
            \draw[->, line width = 0.3mm] (21) -- (22) node[midway, right] {$  |\mathcal{O}_{\mathbb{P}^1}(1) \otimes \mathcal{O}_{\mathbb{P}^2}(1) |$ };
        \end{tikzpicture}
    \end{center} 
    We claim that this diagram commutes. 
    This follows from seeing that 
    \begin{equation}\begin{aligned}
        H^0 (G) \otimes H^0(H) & \rightarrow H^0(G+H) \\
         s \otimes t & \mapsto st
    \end{aligned}\end{equation}
    is a bijection. 
    To prove the claim it is sufficient to show that the kernel is trivial since both domain and codomain have dimension 6. 
    Let $s_0, s_1 \in H^0(G) $ be a basis, and let $S_i = s_i ^{-1}(0)$ be two fibres of the free pencil. 
    Suppose that $s_0 \otimes t_0  + s_1 \otimes t_1 $ belongs to te kernel. 
    Firstly, $t_i$ must both be non-zero. 
    If $t_0 = 0$, then $\forall x \in \Sigma \setminus S_1 $, $t_1 (x) = 0$ which is open dense and so $t_1 =0$ contradicting $r \neq 0$. 
    Analogous contradiction is obtained by taking $t_1 = 0$. 

    On $S_0$, $s_1$ is non-zero, so $t_1|_{S_0} = 0 $. 
    Thus $S_0 \subset t_1 ^{-1}(0)$. 
    Hence $H-G$ is also effective.  
    As $H$ is nef, we arrive then at the contradiction $H. (H-G) =-1$. 
\end{proof}

\subsection{Illustrative examples} 
Table \ref{table_examples} contains the data in the construction of the two pairs of TCS manifolds for which $\xi$ distinguishes the homotopy class of their $\lieGroup{G}_2$-structures. 
The columns consist of: the invariants $b_3$, $m$, and $\xi$; the semi-Fanos; the generators $A_\pm$ of $P_\pm$ respectively; the quadratic form on $P$; and the lattices $\Lambda_\pm$. 
\begin{table}[H]
\caption{Data of the configurations} 
\begin{center} 
\begin{tabular}{cccccccccc}
    $ b_3$ & $ m $ & $ \xi $ &  $ \#_+ $   &  $ \#_- $   &  $ A_+ $  &  $ A_- $  & $P$ & $\Lambda_+$ & $\Lambda_-$\\\hline\\[-5pt] 
            $ 71$ & $ 6 $ & $ 0 $  &   \ref{example_prebb_1_10}   &   \ref{example_prebb_1_10}  &  $ \left( \begin{smallmatrix}1 \end{smallmatrix} \right)$  &   $ \left( \begin{smallmatrix}1 \end{smallmatrix} \right)$  &   $ \left( \begin{smallmatrix} 18 & 0 \\0 & 18 \end{smallmatrix} \right)$  &  $ \left( \begin{smallmatrix}18 \end{smallmatrix} \right)$  &   $ \left( \begin{smallmatrix}18 \end{smallmatrix} \right)$    \\[10pt]
                &  &  $ 24 $  &   \ref{example_prebb_1_Q_0_6} &   \ref{example_prebb_1_PP3_11_14} &   $ \left( \begin{smallmatrix}1 & 1\end{smallmatrix} \right)$  &   $ \left( \begin{smallmatrix}1 & 1\end{smallmatrix} \right)$  &   $ \left( \begin{smallmatrix}16 & 12 & 1 & -1\\12 & 6 & -1 & 1\\1 & -1 & 2 & 5\\-1 & 1 & 5 & 4\end{smallmatrix} \right)$  &   $ \left( \begin{smallmatrix}16 & 12 & 16\\12 & 6 & -16\\16 & -16 & -272\end{smallmatrix} \right)$  &   $ \left( \begin{smallmatrix}2 & 5 & 23\\5 & 4 & -23\\23 & -23 & -552\end{smallmatrix} \right)$ \\[15pt] \hline \\[0pt]
                    $ 85$ & $ 24 $ & $ 12 $  &   \ref{example_prebb_1_PP3_11_14} &   \ref{example_prebb_1_PP3_11_14} &   $ \left( \begin{smallmatrix}1 & 1\end{smallmatrix} \right)$  &   $ \left( \begin{smallmatrix}1 & 1\end{smallmatrix} \right)$  &   $ \left( \begin{smallmatrix}2 & 5 & 1 & -1\\5 & 4 & -1 & 1\\1 & -1 & 2 & 5\\-1 & 1 & 5 & 4\end{smallmatrix} \right)$  &   $ \left( \begin{smallmatrix}2 & 5 & 16\\5 & 4 & -16\\16 & -16 & -272\end{smallmatrix} \right)$  &   $ \left( \begin{smallmatrix}2 & 5 & 16\\5 & 4 & -16\\16 & -16 & -272\end{smallmatrix} \right)$   \\[10pt]
                        & &  $ 36 $  &   \ref{example_prebb_1_PP1_x_PP2} &   \ref{example_prebb_1_PP3_8_2} &   $ \left( \begin{smallmatrix}1 & 2\end{smallmatrix} \right)$  &   $ \left( \begin{smallmatrix}1 & 1\end{smallmatrix} \right)$  &   $ \left( \begin{smallmatrix}0 & 3 & 2 & -2\\3 & 2 & -1 & 1\\2 & -1 & 2 & 8\\-2 & 1 & 8 & 4\end{smallmatrix} \right)$  &   $ \left( \begin{smallmatrix}0 & 3 & 22\\3 & 2 & -11\\22 & -11 & -308\end{smallmatrix} \right)$  &   $ \left( \begin{smallmatrix}2 & 8 & 20\\8 & 4 & -20\\20 & -20 & -180\end{smallmatrix} \right)$   \\[10pt]
\hline
\end{tabular}
\end{center} 
    \label{table_examples} 
\end{table} 

The anticanonical models of \ref{example_prebb_1_PP1_x_PP2} and \ref{example_prebb_1_10} are themselves since they are Fano. 
With a little more care in the other three cases, one establishes that none of the semi-Fanos have anticanonical models
appearing on the list \cite[Theorem 1.1]{jahnke2006gorenstein}. 
Thus they all lead to building blocks of semi-Fano type. 
For each building block, $\mathrm{ker}(H^2(Z) \rightarrow H^2(\Sigma))$ is generated by $c_1(Z)$.  
For each TCS $M$, $N_0 = 0$ so $b_2(M) = 0$.
The cohomology of the building blocks that appear are torsion free while in the arrangements above $P$, $N_\pm + T_\mp$ are all primitive in $L$. 
By Equation \ref{eq_torsion_part_of_H_M}, we conclude that the cohomology of $M$ is also torsion free. 
Hence $(b_3, m, \mu)$ is a complete invariant of the spin diffeomorphism class of $M$, and $\xi(M) \in \mathbb{Z} / (3\tilde{m}) \mathbb{Z}$. 
In each example, $\mu$ can be determined by $\xi$ and agrees for members of a pair, hence it has been omitted. 
It remains then just to check the hypotheses of the genericity results for the lattices $\Lambda_\pm$.
By exhaustive calculation we find the hypotheses are met.

\bibliography{../mybib}
\bibliographystyle{../amsinitial}

\end{document}